\documentclass[12pt]{amsart}
\usepackage{fullpage}
\usepackage{latexsym}
\usepackage{amssymb}

\input amssym.def
\input amssym.tex

\def\page#1{\leaders\hbox to 5mm{.}\hfill \rlap{\hbox to 2mm{\hfill #1}}\par}

\font\tengoth=eufm10 at 12pt
\font\sevengoth=eufm7
\font\fivegoth=eufm5
\newfam\gothfam
\textfont\gothfam=\tengoth
\scriptfont\gothfam=\sevengoth
\scriptscriptfont\gothfam=\fivegoth

\font\cal=cmsy9 at 12pt

\def\og{\leavevmode\raise.3ex\hbox{$\scriptscriptstyle\langle\!\langle\,$}}
\def\fg{\leavevmode\raise.3ex\hbox{$\scriptscriptstyle\,\rangle\!\rangle\ $}}

\newtheorem{theorem}{Theorem}
\newtheorem{cor}[theorem]{Corollary}
\newtheorem{lemma}[theorem]{Lemma}
\newtheorem{prop}[theorem]{Proposition}
\newtheorem{quest}[theorem]{Question}

\newtheorem{definition}[theorem]{Definition}

\newtheorem{fact}[theorem]{Fact}
\newtheorem{claim}[theorem]{Claim}
\newtheorem{notation}[theorem]{Notation}

\newtheorem{assumption}[theorem]{Assumption}
\newtheorem{rk}[theorem]{Remark}
\newtheorem{rks}[theorem]{Remarks}
\newtheorem{example}[theorem]{Example}

\newtheorem{pict}[theorem]{Picture}

\begin{document}

\title{Commensurators of some non-uniform tree lattices and Moufang twin trees}
\author{Peter Abramenko and Bertrand R\'emy}
\maketitle

\vskip 10mm

\centerline{Pr\'epublication de l'Institut Fourier n$^0$ {\bf 627} (2003)}
\centerline{\tt http:$/\!\!/$www-fourier.ujf-grenoble.fr/prepublications.html}

\vskip 10mm

{\footnotesize
{\bf Abstract:~}
Sh. Mozes showed that the commensurator of the lattice
${\rm PSL}_2 \bigl( {\bf F}_p[t{}^{-1}] \bigr)$ is dense in the full automorphism
group of the Bruhat-Tits tree of valency $p+1$, the latter group being much bigger than
${\rm PSL}_2 \bigl({\bf F}_p(\!(t)\!) \bigr)$.
By G.A. Margulis' criterion, this density is a generalized arithmeticity result.
We show that the density of the commensurator holds for many tree-lattices among those called
of Nagao type by H. Bass and A. Lubotzky.
The result covers many lattices obtained via Moufang twin trees.


}

\section*{Introduction}
\label{s - intro}

Tree lattices have been and still are the subject of a lot of interesting
mathematical research, see in particular the monograph \cite{BL} and the papers
cited there.
One of the important issues is to find analogies -- or to establish the differences --
between lattices in the automorphism groups of locally finite trees and lattices in
semisimple Lie groups (i.e. in semisimple algebraic groups over local fields).
These questions are particularly relevant when the Lie groups are of rank 1.
A fundamental theorem due to G.A. Margulis characterizes arithmetic lattices among
lattices in semisimple Lie groups as those which are of infinite index in their
commensurators \cite[Chapter IX, Theorem B]{Margulis}.
(G.A. Margulis proved this characterization for finitely generated lattices, which
excluded non cocompact lattices in rank 1 semisimple Lie groups over local fields
with positive characteristic.
A proof for this remaining case was later given by L. Lifschitz \cite{Lif}.)
G.A. Margulis also showed that being of infinite index in its commensurator
already means that the commensurator of the lattice is \og essentially\fg
dense in the semisimple Lie group, and is in fact dense if the Lie group is simply
connected (for a precise statement, see \cite[Chapter IX, Lemma 2.7]{Margulis}).

\vskip 3mm

These results in the Lie group case motivated (among other things) the study of
commensurators of tree lattices.
In order to fix the ideas, we shall introduce some notation now.
Let $T$ be a locally finite tree with vertex set $VT$, and let $G:={\rm Aut}(T)$ be its
group of automorphisms.
Provided with the usual topology (a basis of open neighborhoods of the identity are the
fixators of finite subtrees of $T$), the group $G$ is locally compact.
A {\it $T$-lattice~} is by definition a lattice in $G$, i.e. a discrete subgroup of
$G$ of finite covolume.
It is standard knowledge \cite[Sections 1.5 and 3.2]{BL} that a subgroup $\Gamma$
of $G$ is a lattice in $G$ if and only if all stabilizers $\Gamma_x$, $x \!\in\! VT$, are
finite, as well as the sum
${\rm Vol}(\Gamma \setminus \!\! \setminus T) := \sum_{x\in\Gamma \setminus
VT}{\mid\!\Gamma_x\!\mid^{-1}}$ is finite.
We call $\Gamma$ a  {\it uniform~} (resp. {\it non-uniform}) $T$-lattice if the
quotient $\Gamma \setminus T$ is finite (resp. infinite).
The {\it commensurator} of $\Gamma$ in $G$ is the group defined as:
${\rm Comm}_G(\Gamma) := \{g \!\in\! G \mid \Gamma \cap g \Gamma g^{-1}$ has finite index in
$\Gamma$ and in $g\Gamma g^{-1} \}$ and will be abbreviated by $C(\Gamma)$ in the
following. Interesting questions concerning $C(\Gamma)$ are the following:

\begin{quest}
\label{q - infinite index?}
Is $\Gamma$ of infinite index in $C(\Gamma)$?
\end{quest}

\begin{quest}
\label{q - dense?}
Is $C(\Gamma)$ (essentially) dense in $G$?
\end{quest}

As mentioned above, these two questions are equivalent for lattices in semisimple
Lie groups, but it is known that they are not equivalent for tree lattices
\cite[Section 10.3]{BL}.
So it is not quite clear which of the two conditions should be used to define, by
analogy, \og arithmetic\fg tree lattices.
Usually one chooses the stronger condition and considers a tree lattice $\Gamma$ as
{\it arithmetic~} if $C(\Gamma)$ is dense in $G$.
An important theorem proved by Y. Liu states that all uniform tree lattices are
arithmetic in the latter sense \cite{Liu}.
Much less is known about commensurators of non-uniform tree lattices.
In particular, to the best of our  knowledge, only two examples of non-uniform tree
lattices with dense commensurators  are discussed in the literature, namely the
example given in \cite[Section 8.3]{BurMoz} and the Nagao lattice
${\rm PSL}_2({\bf F}_p[t^{-1}])$ which is shown in \cite{Moz} to have
a dense commensurator for any prime number $p$.

\vskip 3mm

It is the main objective of the present paper to generalize the two last mentioned
examples in two directions. Firstly, both examples are {\it lattices of Nagao type~}
in the sense of \cite[Chapter 10]{BL}.
For a tree lattice $\Gamma$ of Nagao type, we have a natural {\it level function~}
$\ell:VT \rightarrow {\bf N}$ which is
$\Gamma$-invariant (see Definition \ref{def - level} below), and we set
$L := \{g \!\in\! G \mid \ell(g.x) = \ell(x)$ for all $x \!\in\! VT\}$.
Now the following question is crucial with respect to the examples discussed in
\cite{BurMoz} and \cite{Moz}.

\begin{quest}
\label{q - dense in level?}
Is $C(\Gamma) \cap L$ dense in $L$?
\end{quest}

In Section \ref{s - proof} of this paper, we shall give a positive answer to this question
for all  lattices of {\it directly split Nagao type}, a class of lattices which we shall
introduce in Section \ref{s - Nagao} below.
This is our:

\begin{theorem}
\label{th - density}
If $\Gamma$ is a tree lattice of directly split Nagao type, then $C(\Gamma) \cap L$ is
dense in $L$.
\end{theorem}

Referring the reader to Section \ref{s - Nagao} for some technical details concerning
lattices of directly split Nagao type, we just mention that this class of lattices is
much larger than the examples discussed in \cite{BurMoz} and \cite{Moz}.
In particular, we allow arbitrary finite \og root groups\fg whereas the root groups in
[loc. cit.] are always cyclic.
(For the classical Nagao lattice ${\rm PSL}_2({\bf F}_q[t^{-1}])$ the root groups are
isomorphic to the additive group of ${\bf F}_q$, hence cyclic if and only if $q$ is a prime
number.)

\vskip 3mm

Now if the tree $T$ is not biregular, then one easily checks that $L = G$ (Lemma
\ref{lemma - L equals G}), and hence a positive answer to Question \ref{q - dense in
level?} already means that $C(\Gamma)$ is dense  in $G$.
If, however, $T$ is biregular (as in the case of the classical Nagao lattice), then
we need a second ingredient in order to deduce the density of $C(\Gamma)$
from Theorem \ref{th - density}.
Generalizing a strategy used by Sh. Mozes in \cite{Moz}, we can show that $C(\Gamma)$
is already essentially dense in $G$ in the biregular case if it is not contained in
$L$. More precisely, denoting by $G^{\circ}$ the subgroup (of index at most 2) of $G$
of all type--preserving automorphism of $T$ (preserving the 2--colouring of $T$), we
obtain:

\begin{theorem}
\label{th - biregular density}
If $\Gamma$ is a $T$--lattice of directly split Nagao type with biregular $T$ and
$C(\Gamma)$ is not contained in $L$, then the closure of $C(\Gamma)$ in $G$ contains
$G^{\circ}$.
\end{theorem}

In general it is of course difficult to decide whether $C(\Gamma)$ is included in $L$ or
not. In the case of $\Gamma = {\rm PSL}_2({\bf F}_q[t^{-1}])$, the additional ingredient 
that Sh. Mozes used was the transitive action of $C(\Gamma)$ on the set
of edges of $T$. This transitivity property follows here, for the truly arithmetic group
$\Gamma$, from the fact that its commensurator contains 
${\rm PSL}_2 \bigl( {\bf F}_q(t) \bigr)$, an argument which is not available for other
lattices of Nagao type. This brings us to our second line of generalization.

\vskip 3mm

A special class of lattices of Nagao type arises from the theory of twin trees as follows.
Let a locally finite (thick) twin tree $(T_\pm,\delta^*)$ be given which has the Moufang
property  (for definitions, see Section \ref{s - mtt} below).
We remark that the trees $T_+$ and $T_-$ are necessarily biregular.
Denote by $A$ the automorphism group of the twin tree $(T_\pm,\delta^*)$ and by $\Lambda$ its
subgroup generated by all root groups.
Fix a vertex of $T_-$ and denote by $\Gamma$ its stabilizer in $\Lambda$.
Then it is well-known that $\Gamma$ acts on $T_+$ as a non-uniform $T_+$-lattice with a ray 
$R$ as fundamental  domain.
Moreover, $\Gamma$ is always a lattice of Nagao type in the sense of \cite{BL}, and it is
of directly split Nagao type in our sense if it satisfies the following  condition (which
will also be explained in more detail in Section \ref{s - mtt}).

\begin{assumption}
\label{ass - comm}
{\rm (Comm)~} For any prenilpotent pair $\{a;b \}$ of twin roots of the Moufang
twin tree  $(T_\pm,\delta^*)$, the corresponding root groups $U_a$ and $U_b$ commute.
\end{assumption}

It is not difficult to check that $A$ is contained in the commensurator of the
$T_+$-lattice $\Gamma$ and that it acts transitively on the set of all (geometric) edges
of $T$. In particular, $C(\Gamma)$ is not contained in $L$, and thus Theorem 
\ref{th - biregular density} implies (where $G^{\circ}$ denotes the group of all 
type--preserving automorphisms of $T_+$):

\begin{theorem}
\label{th - dense when Moufang}
If the thick locally finite Moufang twin tree $(T_\pm,\delta^*)$ satisfies {\rm (Comm)~} and
$\Gamma$ is the $T_+$-lattice described above, then the closure of $C(\Gamma)$ contains
$G^\circ$. Therefore $C(\Gamma)$ is dense in $G$ if $T_+$ is not regular and also, in the
regular case, if $A$ contains an automorphism interchanging the types of vertices.
\end{theorem}

First examples for groups $\Lambda$ and $\Gamma$ satisfying (Comm) come from rank 2 Kac-Moody
groups over finite fields (see \ref{ss - examples}).
As a special (affine) case we obtain
$\Lambda = {\rm PSL}_2({\bf F}_q[t,t^{-1}])$ and $\Gamma = {\rm PSL}_2({\bf F}_q[t^{-1}])$
for any prime power $q$.
In the latter case, the matrix
$\left( \begin{array}{rr}t & 0\\ 0 & 1\\ \end{array} \right) \in {\rm PGL}_2({\bf F}_q[t])$
provides an element of $A$ interchanging the vertex types: it acts on the
tree as a hyperbolic translation, and its translation length is the length of a single edge.
The construction of more exotic examples is sketched in Example \ref{ex - not linear but
arithmetic}.
As a very special case, we consider trees of constant degree equal to 7.
For a non-uniform tree-lattice $\Gamma$ naturally defined in terms of twinnings, the density
of the commensurator holds, at least in the index 2 subgroup $G^\circ$.
Note that in this case the full automorphism group $G$ is \og as big as\fg in the case
considered by Sh. Mozes, in the sense that the tree of valency 7 is \og as homogeneous as\fg a
Bruhat-Tits tree of some ${\rm PSL}_2 \bigl({\bf F}_p(\!(t)\!) \bigr)$.
Nevertheless, the case is new at least because it can easily be proved that the lattice
itself cannot be linear over a field.

\vskip 3mm

We close this introduction by an open question relevant to a higher-dimensional
generalization of the above results.

\begin{quest}
Let $\Delta_\pm$ be a locally finite thick Moufang twin building whose apartments are
right-angled tilings of the hyperbolic plane.
Does the non-uniform lattice naturally defined as a negative chamber stabilizer have a dense
commensurator in ${\rm Aut}(\Delta_+)$?
\end{quest}

We focus on the case of right-angled Fuchsian buildings because the local product
structure of the links at vertices implies the existence of many symmetries for this class
of buildings.
Therefore the corresponding full automorphism groups are big topological groups.
Recent results by F. Haglund on commensurators of some uniform lattices for hyperbolic
buildings \cite{HagComm} are analogous to Y. Liu's theorem for trees (and might be used as
such).
Note that concrete examples of right-angled Fuchsian buildings admitting a twinning
are available \cite[Theorem 4.E.2]{RemRon}.

\vskip 3mm

This paper is organized as follows.
Section \ref{s - Nagao} deals with lattices of Nagao type as introduced by H. Bass and A.
Lubotzky.
The subclass of tree-lattices of directly split Nagao type is introduced
(\ref{ss - directly split Nagao type}) and
the condition defining it is characterized in terms of actions on some horoballs
(\ref{ss - horoballs and horospheres}).
Section \ref{s - proof} proves the main density result
(Theorem \ref{th - density}), following the ideas presented in \cite[Section 8.3]{BurMoz}.
The main change with respect to \cite{BurMoz} and \cite{Moz} is the replacement of local data
by purely group-theoretic considerations (\ref{ss - permuting elements of the same level}).
Section \ref{s - mtt} recalls some facts about Moufang twin trees, and then describes the
intersection of this theory with that of tree-lattices of Nagao type.
This shows in particular that the class of non-uniform tree-lattices to which
Theorem \ref{th - biregular density} applies is quite wide (\ref{ss - examples}).

\vskip 3mm

We thank Sh. Mozes for helpful discussions on Example \ref{ex - not linear but arithmetic}
and are very grateful to the organizers of the conference \og Geometric Group Theory\fg 
(Guwahati, Assam -- India, December 2002)  organized by the Indian Institute of Technology 
Guwahati and the Indian Statistical Institute (supported by the National Board of Higher 
Mathematics).
The conference gave the authors the opportunity to complete their research on the subject
which is presented in this paper.

\vskip 6mm

\section{Lattices of Nagao type}
\label{s - Nagao}

The tree lattices we are discussing in this paper are always fundamental groups of
rays of groups.
More precisely, we shall work with the following set-up.

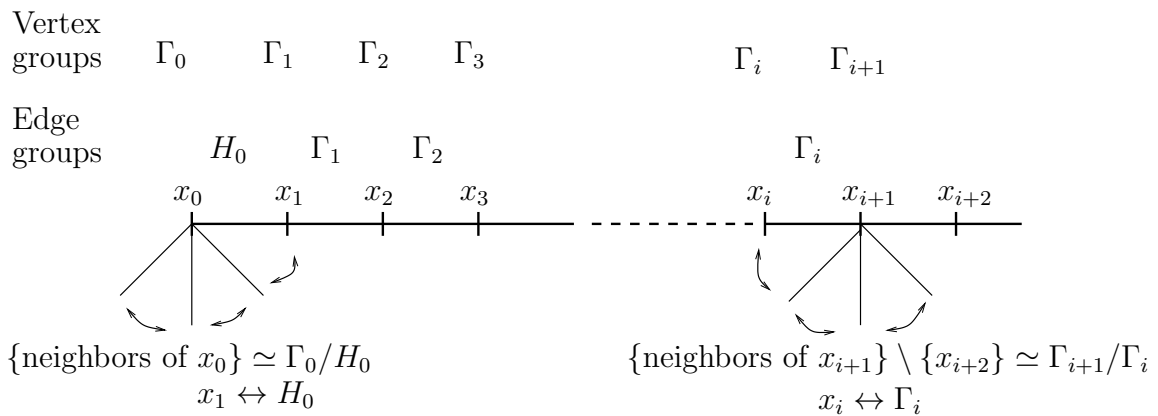
\begin{pict}
\label{pict - graph of groups}
\end{pict}

\vskip 3mm
\centerline{\begin{picture}(0,0)%
\includegraphics{Picture1.pstex}%
\end{picture}%
\setlength{\unitlength}{1973sp}%
\begingroup\makeatletter\ifx\SetFigFont\undefined
\def\x#1#2#3#4#5#6#7\relax{\def\x{#1#2#3#4#5#6}}%
\expandafter\x\fmtname xxxxxx\relax \def\y{splain}%
\ifx\x\y   
\gdef\SetFigFont#1#2#3{%
  \ifnum #1<17\tiny\else \ifnum #1<20\small\else
  \ifnum #1<24\normalsize\else \ifnum #1<29\large\else
  \ifnum #1<34\Large\else \ifnum #1<41\LARGE\else
     \huge\fi\fi\fi\fi\fi\fi
  \csname #3\endcsname}%
\else
\gdef\SetFigFont#1#2#3{\begingroup
  \count@#1\relax \ifnum 25<\count@\count@25\fi
  \def\x{\endgroup\@setsize\SetFigFont{#2pt}}%
  \expandafter\x
    \csname \romannumeral\the\count@ pt\expandafter\endcsname
    \csname @\romannumeral\the\count@ pt\endcsname
  \csname #3\endcsname}%
\fi
\fi\endgroup
\begin{picture}(12783,5088)(76,-6394)
\put(151,-2761){\makebox(0,0)[lb]{\smash{\SetFigFont{12}{14.4}{rm}Edge }}}
\put(151,-3181){\makebox(0,0)[lb]{\smash{\SetFigFont{12}{14.4}{rm}groups }}}
\put(2176,-3661){\makebox(0,0)[lb]{\smash{\SetFigFont{12}{14.4}{rm}$x_0$}}}
\put(3451,-3661){\makebox(0,0)[lb]{\smash{\SetFigFont{12}{14.4}{rm}$x_1$}}}
\put(4576,-3661){\makebox(0,0)[lb]{\smash{\SetFigFont{12}{14.4}{rm}$x_2$}}}
\put(5776,-3661){\makebox(0,0)[lb]{\smash{\SetFigFont{12}{14.4}{rm}$x_3$}}}
\put(9376,-3661){\makebox(0,0)[lb]{\smash{\SetFigFont{12}{14.4}{rm}$x_i$}}}
\put(10576,-3661){\makebox(0,0)[lb]{\smash{\SetFigFont{12}{14.4}{rm}$x_{i+1}$}}}
\put(11776,-3661){\makebox(0,0)[lb]{\smash{\SetFigFont{12}{14.4}{rm}$x_{i+2}$}}}
\put(1951,-1936){\makebox(0,0)[lb]{\smash{\SetFigFont{12}{14.4}{rm}$\Gamma_0$}}}
\put(3301,-1936){\makebox(0,0)[lb]{\smash{\SetFigFont{12}{14.4}{rm}$\Gamma_1$}}}
\put(4501,-1936){\makebox(0,0)[lb]{\smash{\SetFigFont{12}{14.4}{rm}$\Gamma_2$}}}
\put(5701,-1936){\makebox(0,0)[lb]{\smash{\SetFigFont{12}{14.4}{rm}$\Gamma_3$}}}
\put(9226,-2011){\makebox(0,0)[lb]{\smash{\SetFigFont{12}{14.4}{rm}$\Gamma_i$}}}
\put(9976,-3136){\makebox(0,0)[lb]{\smash{\SetFigFont{12}{14.4}{rm}$\Gamma_i$}}}
\put(2626,-3136){\makebox(0,0)[lb]{\smash{\SetFigFont{12}{14.4}{rm}$H_0$}}}
\put(3901,-3136){\makebox(0,0)[lb]{\smash{\SetFigFont{12}{14.4}{rm}$\Gamma_1$}}}
\put(5176,-3136){\makebox(0,0)[lb]{\smash{\SetFigFont{12}{14.4}{rm}$\Gamma_2$}}}
\put(10426,-2011){\makebox(0,0)[lb]{\smash{\SetFigFont{12}{14.4}{rm}$\Gamma_{i+1}$}}}
\put( 76,-5761){\makebox(0,0)[lb]{\smash{\SetFigFont{12}{14.4}{rm}$\{$neighbors of $x_0\} \simeq \Gamma_0/H_0$}}}
\put(151,-1981){\makebox(0,0)[lb]{\smash{\SetFigFont{12}{14.4}{rm}groups}}}
\put(151,-1561){\makebox(0,0)[lb]{\smash{\SetFigFont{12}{14.4}{rm}Vertex }}}
\put(2476,-6211){\makebox(0,0)[lb]{\smash{\SetFigFont{12}{14.4}{rm}$x_1 \leftrightarrow H_0$}}}
\put(10276,-6286){\makebox(0,0)[lb]{\smash{\SetFigFont{12}{14.4}{rm}$x_i \leftrightarrow \Gamma_i$}}}
\put(7876,-5761){\makebox(0,0)[lb]{\smash{\SetFigFont{12}{14.4}{rm}$\{$neighbors of $x_{i+1}\}\setminus \{x_{i+2}\} \simeq \Gamma_{i+1}/\Gamma_i$}}}
\end{picture}
}
\vskip 3mm

\begin{notation}\rm
\label{notation - Nagao}
Let $R$ be an infinite ray with vertices $\{x_i \}_{i \geq 0}$.
To each vertex $x_i, \, i \geq 0$, we attach a finite group $\Gamma_i$ such that
$\Gamma_i$ is a subgroup of $\Gamma_{i+1}$ for all $i > 0$.
Furthermore, we are given a
common subgroup $H_0$ of $\Gamma_0$ and $\Gamma_1$. We  define the integers $k :=
[\Gamma_0 : H_0]$, $q_0 := k - 1$, $q_1 := [\Gamma_1 : H_0]$  and $q_i := [\Gamma_i :
\Gamma_{i-1}]$ for $i \geq 2$, and we require that $q_i \geq 2$  for all $i \geq 0$. Now
we attach $H_0$ to the edge $\{x_0,x_1\}$ of $R$ and $\Gamma_i$  to $\{x_i,x_{i+1}\}$ for
all $i > 0$. We have thus defined a graph of groups, where the  graph is the ray $R$.
\end{notation}

\subsection{Levels, degrees and a first density result}
\label{ss - level and degree}
Denote by $\Gamma$ the fundamental group of this graph of groups and by $T$ the
canonical tree provided by Bass-Serre theory \cite[I.5.1]{Ser}: $\Gamma$ acts on $T$ with
$R$ as a fundamental domain.
For any $i \geq 0$, we identify $\Gamma_i$ with the vertex stabilizer 
${\rm Stab}_\Gamma(x_i)$ and $H_0 = \Gamma_0 \cap \Gamma_1$ with the edge stabilizer 
${\rm Stab}_\Gamma(\{ x_0,x_1 \} )$.

\begin{definition}
\label{def - level}
Any vertex $x \!\in\! VT$ lies in a $\Gamma$-orbit $\Gamma.x_i$ for precisely one $i
\geq 0$; we call $i$ the {\rm level~} of $x$ and denote it by $\ell(x) = i$.
We set $G := {\rm Aut}(T)$ and
$L := \{g \in G \mid \ell(g.x) = \ell(x)$ for all $x \!\in\! VT\}$.
\end{definition}

The kernel of the canonical homomorphism $\Gamma \rightarrow G$ is a (finite) subgroup
of $H_0$, and by a slight abuse of notation, we shall denote the image of $\Gamma$ in
$G$ again by $\Gamma$.  Therefore $\Gamma$ is a discrete subgroup of $G$, and it is of
course also a subgroup of $L$.
Our assumption $q_i \geq 2$ for all $i \geq 0$ implies that
Vol$(\Gamma \setminus \!\! \setminus T) = \sum_{i = 0}^{\infty}|\Gamma_i|^{-1}$ is finite.
Hence $\Gamma$ is an (obviously non-uniform) $T$-lattice, called a {\it lattice of Nagao
type~} in \cite[Chapter 10]{BL}.
Since $R$ is a fundamental domain for the action of $\Gamma$ on $T$, we immediately
obtain the following two facts.

\begin{fact}
\label{fact - valency at level 0}
Each vertex $x \!\in\! VT$ of level $0$ has precisely $k=q_0+1$ neighbors,
and all of them are of level $1$.
In particular, the degree (or valency) of $x$ is ${\rm deg}(x)=k$.
\qed\end{fact}

\begin{fact}
\label{fact - valency at level > 0}
Each vertex $x \!\in\! VT$ with $\ell(x)=i>0$ has precisely $1$ neighbor of level $i+1$ and
$q_i$ neighbors of level $i-1$.
In particular, ${\rm deg}(x)=q_i+1$.
\qed\end{fact}

Fact \ref{fact - valency at level > 0} immediately implies the following two statements.

\begin{fact}
\label{fact - level decreases}
If $\{ y_j \}_{j \geq 0} \subset VT$ is a geodesic ray in $T$ with
$\ell(y_0) > \ell(y_1)$, then for any $j \leq \ell(y_0)$, we have:
$\ell(y_j)=\ell(y_0)-j$.
\qed\end{fact}

\begin{fact}
\label{fact - level increases}
If $y \!\in\! VT$ has level $i > 0$, then there is a unique (infinite) ray
$R(y) = (y_0,y_1, \cdots )$ starting in $y = y_0$ such that $\ell(y_j) = i + j$ for
all $j \geq 0$.
The ray $R(y)$ will be referred to as the {\rm level-increasing ray~} from $x$.
\qed\end{fact}

It is well-known that either $T$ is biregular (i.e. deg($x$) = deg($y$) for any two
vertices $x,y \!\in\! VT$ at even distance from each other) or $G \setminus T \, = \, R$
\cite[Section 10.1]{BL}.
However, we shall give a short elementary argument for this (in fact for a slightly more
general) statement here, also because we shall later need the remark following this lemma.

\begin{lemma}
\label{lemma - L equals G}
If $F$ is a subgroup of $G$ which contains $\Gamma$ and is not included in $L$,
then $F.x_0 = F.x_{2m}$ and $F.x_1 = F.x_{2m+1}$ for all integers $m \geq 0$.
Hence in this case $T$ has to be biregular, and $F$ acts transitively on vertices of
the same type (i.e. at even distance from each other) and on the set $ET$ of
(geometric) edges of $T$.
\end{lemma}

{\it Proof.~} Since for any $x \!\in\! VT$, $\ell(x)$ is equal to the distance of $x$
from the orbit $\Gamma.x_0$, $F$ acts level-preservingly on $T$ if it stabilizes
$\Gamma.x_0$.
Since $F$ is not included in $L$, there is an $n > 0$ with $F.x_0 = F.x_n$. We
choose $n > 0$ minimal with this property. If $n = 1$, we inductively get
$F.x_r = F.x_{r-1} = F.x_0$ for all $r \!\in\! {\bf N}$. If $n = 2$, then inductively
$F.x_3 = F.x_1$, $F.x_4 = F.x_0 (= F.x_2)$, $F.x_5 = F.x_1$ and so on.

\vskip 3mm

Now assume that $n > 2$. Since $F.x_{n-1} = F.x_1$, we can choose a $g \!\in\! F$ with
$g.x_{n-1} = x_1$. Since $g.x_n \neq x_2$ (otherwise $F.x_2 = F.x_n = F.x_0$,
contradicting the minimality of $n$), $\ell(g.x_n) = 0$
by Fact \ref{fact - valency at level > 0}. Now $x_1$ has
$q_1 \geq 2$ neighbors of level 0, and $\Gamma_1$ acts transitively on them. Hence we
can choose $\gamma_1 \!\in\! \Gamma_1$ such that $\gamma_1.(g.x_n) \neq g.x_n$. Hence
$f := g^{-1}\gamma_1g$ is in $F_{x_{n-1}}$ and $f.x_n \neq x_n$. Therefore
$\Gamma_{n-1}f.x_n$ contains $x_{n-2}$, and thus $F.x_{n-2} = F.x_n = F.x_0$.
However, this again contradicts the minimality of $n$, and so $n > 2$ is impossible.

\vskip 3mm

This proves the first claim, which immediately shows that $F$ acts transitively on
vertices of the same type. In particular, for any edge $e \!\in\! ET$, there is an
$f \!\in\! F$ such that $x_0 \in f.e$, and then there is a $\gamma_0 \!\in\! \Gamma_0$
such that $(\gamma_0.f).e = \{ x_0,x_1 \} $.
Hence $F$ acts transitively on $ET$.
\qed

\begin{rk}\rm
\label{rk - biregular}
It is not needed in the above argument that $R$ is infinite.
So if we have the same set-up as above, only with a finite path $R'$ instead of $R$, and
if $T'$ and $G'$ are defined similarly, then also $T'$ has to be biregular or else
$G' \! \setminus \! T' = R'$.
\end{rk}

\begin{lemma}
\label{lemma - extension to L}
Any level-preserving isomorphism $\psi: T_1 \rightarrow T_2$ between two
subtrees $T_1$ and $T_2$ of $T$ extends to an element of $L$.
\end{lemma}

{\it Proof.~}
By Zorn's lemma it suffices to show that $\psi$ can always be properly and
level-preservingly extended if $T_1 \neq T$.
That this is in fact true can  readily be checked by using the Facts
\ref{fact - valency at level 0} and \ref{fact - valency at level > 0}
above.
\qed

\begin{rk}\rm
This lemma in particular shows that the group $L$ is non-discrete and uncountable.
\end{rk}

The elementary statements which we deduced above have an interesting consequence if
we combine them with an idea presented in the last paragraph of Mozes' paper \cite{Moz}.
Recall that $G^\circ$ denotes the group of type-preserving automorphisms of $T$.

\begin{prop}
\label{prop - maximality of L}
If $F$ is a subgroup of $G$ properly containing $L$ (so that $T$ must be
biregular), then the closure of $F$ in $G$ contains $G^\circ$.
\end{prop}

{\it Proof.~}
We have to show that any type-preserving isomorphism
$\phi: S_1 \rightarrow S_2$ between two finite subtrees $S_1$ and $S_2$ of $T$ can
be extended to an element of $F$. Without loss of generality, we may assume that
$S_1 = B_s(z_1)$, $S_2 = B_s(z_2)$ are two balls of radius $s \in {\mathbb N}$
with centers $z_1$ and $z_2$ which are of the same type. Let $v_1$ be a terminal
vertex in $S_1$ and set $v_2 = {\phi}(v_1)$. So also $v_1$ and $v_2$ are of the same
type. Let $u_j$ ($j = 1, \, 2$) be the unique neighbor of $v_j$ in $S_j$. Choose a
vertex $x_n \in R$ of the same type as $v_1$ and $v_2$ with $n \geq 2s$.
By Lemma \ref{lemma - L equals G}, there exist elements $f_j \in F$
satisfying $f_j.v_j = x_n$ and $f_j.u_j = x_{n-1}$ for $j = 1, \, 2$.
Set $T_j := f_j(S_j)$ for $j = 1, \, 2$, and consider the isomorphism
$\psi:= f_2 \! \mid_{S_2}{\phi} {f_1}^{-1} \! \mid_{T_1}: T_1 \rightarrow T_2$.
Note that $\psi$ fixes $x_n$ and $x_{n-1}$ by construction.
So Fact \ref{fact - level decreases} together with
$n \geq 2s$ implies that $\psi$ is automatically level-preserving. Hence by
Lemma \ref{lemma - extension to L}, there exists a $g \in L$ satisfying
$g \! \mid_{T_1}  = \psi$.
Therefore the element $f:= {f_2}^{-1}gf_1 \in F$ extends $\phi$.
\qed

\subsection{Lattices of directly split Nagao type}
\label{ss - directly split Nagao type}
Recall that we want to investigate the commensurator
$C(\Gamma)=\{g \in G \mid \Gamma \cap g \Gamma g^{-1}$ has finite index in $\Gamma$
and in $g\Gamma g^{-1} \}$.
In the present generality, we even cannot exclude the possibility that $C(\Gamma) = \Gamma$
(for such an example, see \cite[Section 10.5]{BL}).
Therefore we have to consider a restricted class of lattices of Nagao type.
Motivated by and generalizing the discussion in \cite[Section 8.3]{BurMoz}, we choose the
following one.

\begin{definition}
\label{def - directly split Nagao}
We say that the above introduced tree lattice $\Gamma$ is of
{\rm directly split Nagao type~} if there exist, for all $j > 0$, $H_0$-invariant subgroups
$U_j \leq \Gamma_j$, called {\rm root groups~}, such that
$\Gamma_i = H_0 \ltimes (U_1 \times U_2 \times \ldots \times U_i)$ for all $i > 0$.
\end{definition}

\begin{rks}\rm
1. \og Directly split\fg above refers to the direct products $U_1 \times
U_2 \times \ldots \times U_i$.
In view of this terminology one would perhaps expect that also the product with $H_0$ is
direct.
However, this would be too restrictive; it would even exclude the classical Nagao lattice
${\rm PSL}_2({\bf F}_q[t^{-1}])$.

2. The terminology \og root group\fg is motivated by building theory (here: twin tree
theory) and will be explained in Section \ref{s - mtt}.

3. If $H_0 = \{1\}$, our lattices of directly split Nagao type are the same as the
{\it product groupings~} associated to Nagao rays in \cite[Section 10.6]{BL}.

4. We stress again that our root groups $U_i$ can be completely arbitrary
non-trivial finite groups.

5. It would be desirable, especially in view of Moufang twin trees which do not
satisfy the Condition (Comm),  to obtain results as deduced below for splittings $\Gamma_i
= H_0 \ltimes (U_1 \cdots U_i)$ which are not necessarily direct. However, our method
really requires that the different root groups commute.
We shall give a geometric characterization of this algebraic property in Lemma
\ref{lemma - action on std hb} below.
\end{rks}

For the rest of this paper, we shall assume that our lattices are of directly split
Nagao type.
We start by collecting a few elementary statements, the first two of which are obvious.

\begin{fact}
\label{fact - size of root gp}
We have:
$\Gamma_1 = H_0 \ltimes U_1$ and $\Gamma_i = \Gamma_{i-1} \ltimes U_i$
for any $i \geq 2$.
In particular, $\mid\!U_i\!\mid = q_i \geq 2$ for any $i \geq 1$.
\qed\end{fact}

\begin{fact}
\label{fact - transitivity of root gp}
For any $i > 0$, the group $U_i$ acts simply transitively (i.e. regularly) on the set
of neighbors of $x_i$ which are different from $x_{i+1}$.
\qed\end{fact}

This can easily be generalized as follows.

\begin{lemma}
\label{lemma - products}
For any two integers $i$ and $j$ with $0 < i \leq j$, the group
$U_{i,j} := U_i \times U_{i+1} \ldots \times U_j$ acts simply transitively on
$M_{i,j} := \{ x \!\in\! VT \mid \ell(x) = i-1$ and ${\rm dist}(x, x_j)=j-i+1 \}$,
where ${\rm dist}(x, x_j)$ denotes the distance between $x$ and $x_j$ in $T$.
\end{lemma}

{\it Proof.~}
The fact that $U_{i,j}$ acts transitively on $M_{i,j}$ immediately follows
from Fact \ref{fact - transitivity of root gp}
by induction on $j-i$. Now assume that an element
$u = u_j{\ldots}u_i \in U_{i,j}$ fixes $x_{i-1}$. Since $U_{i,j} \leq \Gamma_j$,
$u$ has to fix $x_l$ for all $i-1 \, \leq l \, \leq j$. In particular,
$x_{j-1} = u.x_{j-1} = u_j.x_{j-1}$, which implies $u_j = 1$, again by Fact
\ref{fact - transitivity of root gp}.
Going on this way, we get $u_j = u_{j-1} = \ldots = u_i = 1$.
\qed

\subsection{Horoballs and horospheres}
\label{ss - horoballs and horospheres}
Recall that to any end $\epsilon$ of the tree $T$, we can associate horospheres and
horoballs centered at $\epsilon$, see for instance \cite[Section 9.2]{BL}.
In the present paper, only certain horospheres and horoballs, defined by vertices of
positive level, will play a role.

\begin{definition}
Let $x$ be a vertex of level $\ell(x) > 0$, and let $\epsilon_x$ be
the end of $T$ defined by the level-increasing ray from $x$ (see Fact \ref{fact - level
increases}).
\begin{enumerate}
\item[(a)] We denote by ${\rm HS}(x)$ the horosphere centered at $\epsilon_x$
which contains $x$.
\item[(b)] We denote by ${\rm HB}(x)$ the horoball centered at $\epsilon_x$
which has ${\rm HS}(x)$ as its boundary, i.e. as its set of terminal vertices.
\end{enumerate}
\end{definition}

\begin{rk}\rm
It is easy to see that the horoball ${\rm HB}(x)$ is the connected component of $x$ in the
forest spanned by all vertices of level $\geq \ell(x)$, and that the horosphere ${\rm HS}(x)$
is the set of vertices of level $\ell(x)$ in ${\rm HB}(x)$.
\end{rk}

The following statement is an immediate consequence of Lemma \ref{lemma - products}
and the above definition.

\begin{fact}
\label{fact - action on std hs}
For any integer $i > 0$, the subgroup $W_{i+1} := \langle U_j \mid j \geq i + 1 \rangle$
of $\Gamma$ acts simply transitively on the horosphere ${\rm HS}(x_i)$.
\qed\end{fact}

\begin{rk}\rm
As opposed to Lemma \ref{lemma - action on std hb} below,
this fact remains true in a more general context (for instance for general
Moufang twin trees), where the products $U_1{\cdots}U_i$ are not required to be direct.
\end{rk}

We close this section by stating a geometric consequence of the directness of these
products which will play an important role in the next section.

\begin{lemma}
\label{lemma - action on std hb}
For any $i > 0$, the group $U_i$ fixes the horoball ${\rm HB}(x_i)$ pointwise.
\end{lemma}

{\it Proof.~}
It suffices to show that $U_i$ acts trivially on ${\rm HS}(x_i)$. So take an
arbitrary vertex $x \!\in\! {\rm HS}(x_i)$, and let $w$ be the unique element in
$W_{i+1}$ satisfying $x = w.x_i$ (Fact \ref{fact - action on std hs}).
Now since $U_i$ commutes with
$U_j$ for all $j > i$, it also commutes with $W_{i+1}$. Therefore
$u_i.x = u_i.(w.x_i) = w.(u_i.x_i) = w.x_i = x$ for all $u_i \!\in\! U_i$.
\qed\vskip 3mm

\begin{pict}
\label{pict - horoballs}
\end{pict}

\vskip 3mm
\centerline{\begin{picture}(0,0)%
\includegraphics{Picture2.pstex}%
\end{picture}%
\setlength{\unitlength}{1973sp}%
\begingroup\makeatletter\ifx\SetFigFont\undefined
\def\x#1#2#3#4#5#6#7\relax{\def\x{#1#2#3#4#5#6}}%
\expandafter\x\fmtname xxxxxx\relax \def\y{splain}%
\ifx\x\y   
\gdef\SetFigFont#1#2#3{%
  \ifnum #1<17\tiny\else \ifnum #1<20\small\else
  \ifnum #1<24\normalsize\else \ifnum #1<29\large\else
  \ifnum #1<34\Large\else \ifnum #1<41\LARGE\else
     \huge\fi\fi\fi\fi\fi\fi
  \csname #3\endcsname}%
\else
\gdef\SetFigFont#1#2#3{\begingroup
  \count@#1\relax \ifnum 25<\count@\count@25\fi
  \def\x{\endgroup\@setsize\SetFigFont{#2pt}}%
  \expandafter\x
    \csname \romannumeral\the\count@ pt\expandafter\endcsname
    \csname @\romannumeral\the\count@ pt\endcsname
  \csname #3\endcsname}%
\fi
\fi\endgroup
\begin{picture}(11700,6805)(451,-6769)
\put(2176,-4561){\makebox(0,0)[lb]{\smash{\SetFigFont{12}{14.4}{rm}$x_i$}}}
\put(3301,-4561){\makebox(0,0)[lb]{\smash{\SetFigFont{12}{14.4}{rm}$x_{i+1}$}}}
\put(976,-4561){\makebox(0,0)[lb]{\smash{\SetFigFont{12}{14.4}{rm}$x_{i-1}$}}}
\put(451,-3436){\makebox(0,0)[lb]{\smash{\SetFigFont{12}{14.4}{rm}$U_i$}}}
\put(2101,-3286){\makebox(0,0)[lb]{\smash{\SetFigFont{12}{14.4}{rm}$U_{i+1}$}}}
\put(10501,-4561){\makebox(0,0)[lb]{\smash{\SetFigFont{12}{14.4}{rm}$x_{j-1}$}}}
\put(12151,-3961){\makebox(0,0)[lb]{\smash{\SetFigFont{12}{14.4}{rm}$x_j$}}}
\put(6301,-1336){\makebox(0,0)[lb]{\smash{\SetFigFont{12}{14.4}{rm}length $j-i+1$}}}
\end{picture}
}
\vskip 6mm

\section{Proof of the main density theorem}
\label{s - proof}

We are now going to investigate the commensurator $C(\Gamma)$, where $\Gamma$ is a
$T$-lattice of directly split Nagao type as introduced in Section \ref{s - Nagao}.
Let us first of all mention that Question \ref{q - infinite index?} formulated in the
introduction has a positive answer in this case.
In fact, according to \cite[Section 10.1]{BL}, $\Gamma$ is already of infinite index in
$N_G(\Gamma)$, its normalizer in $G$.
In this section, it is our goal to show that $C(\Gamma) \cap L$ is dense in $L$, which is
the major step towards answering Question \ref{q - dense?}. Our strategy will be an
appropriate adaptation (using more group theoretic constructions instead of \og local data\fg)
of the method developed in \cite[Section 8.3]{BurMoz}. In particular, we shall also prove an
extension result for commensurators of certain uniform tree lattices for subtrees of $T$ with
bounded level function (see Proposition \ref{prop - extension-comm}) and then apply Liu's
theorem (\cite{Liu}).

\vskip 3mm

As in \cite{Moz}, it will be technically convenient to work with a certain finite index
subgroup $\Delta$ of $\Gamma$ (so that of course $C(\Delta) = C(\Gamma)$) instead of
$\Gamma$ itself.
In the case of the classical Nagao lattice ${\rm PSL}_2({\bf F}_q[t^{-1}])$, Sh. Mozes used
the first congruence subgroup.
In the theory of Kac-Moody groups, one would analogously choose $\Delta$ to be the
\og unipotent radical\fg of $\Gamma$, where $\Gamma$ is considered as a parabolic
subgroup of a Kac-Moody group over a finite field.
In the next subsection we shall define $\Delta$ in terms of the subgroups $U_i$ introduced
in Definition \ref{def - directly split Nagao}.

\subsection{The group $\Delta$ and its action on the tree and on horoballs}
\label{ss - Delta-action on the tree and on horoballs}
We keep the notations introduced in Section \ref{s - Nagao} and define further subgroups of
the $T$-lattice $\Gamma$ of directly split Nagao type.
First, we call $V:= \langle U_i\mid i>0 \rangle = \bigoplus_{i>0} U_i$ the
direct sum of the groups $U_i$ when $i$ ranges over all the positive integers $i$.
We fix a system $\{\gamma_1:=1$; $\gamma_2$; $\cdots$ ; $\gamma_k\}$
of representatives of the cosets in $\Gamma_0/H_0$.
For $i>0$ and $1 \leq s \leq k$, we also set:

\vskip 3mm

\centerline{$U_{i,s} := \gamma_sU_i\gamma_s{}^{-1}$ \qquad and \qquad
$V_s := \langle U_j,s \mid j>0 \rangle = \gamma_sV\gamma_s{}^{-1}$,}

\vskip 3mm

so that $V_1=V$, and we finally define $\Delta := \langle V_s \mid 1 \leq s \leq k
\rangle$.

\begin{lemma}
\label{lemma - Levi}
The group $\Delta$ is the normal closure of $V$ in $\Gamma$, and we have:
$\Gamma=\Gamma_0 \ltimes \Delta$.
\end{lemma}

{\it Proof.~}
Recall that $\Gamma = \langle \Gamma_0 \cup V \rangle$ and that $H_0$ normalizes
$V$. This implies $\Gamma = \Gamma_0\Delta$ (since $V\Gamma_0\Delta = \Gamma_0\Delta$),
and hence $\Delta$ is normal in $\Gamma$.

\vskip 3mm

The definition of $\Gamma$ as the fundamental group of a ray of groups implies
that $\Gamma$ is an amalgam of the form $\Gamma = \Gamma_0 \, *_{H_0} \, H_0V$.
Now the normal form for amalgams (as stated for instance in \cite[I.1.2]{Ser})
shows that no product of the form
$\gamma_{j_1}v_1\gamma_{j_1}{}^{-1}\gamma_{j_2}v_2\gamma_{j_2}{}^{-1}{\ldots}
\gamma_{j_l}v_l\gamma_{j_l}{}^{-1}$ with $l \! \in \! {\mathbb N}$,
$v_1$, $\ldots$, $v_l \! \in \! V \setminus \{ 1 \}$, $ 1 \leq j_r \leq k$ for all
$r \leq l$ and $j_r \neq j_{r+1}$ (hence $\gamma_{j_r}{}^{-1}\gamma_{j_{r+1}} \!
\not\in  \! H_0$) for all $r < k$ is an element of $\Gamma_0$.
Hence $\Gamma_0 \cap \Delta = \{ 1 \}$.
\qed\vskip 3mm

The same normal form argument also shows that $\Delta = V_1 * V_2 * ... *V_k$.

\begin{lemma}
\label{lemma - action on std hs and hb}
For any positive integer $i > 0$, the stabilizer ${\rm Stab}_{\Delta}(x_i)$ is equal to
$U_1 \times U_2 \times ... \times U_i$ and acts trivially on ${\rm HB}(x_i)$.
\end{lemma}

{\it Proof.~}
Since $H_0 \cap \Delta = \{ 1 \} $ by the previous lemma, we obtain
${\rm Stab}_{\Delta}(x_i) = \Gamma_i \cap \Delta =
\bigl( H_0 \ltimes (U_1 \times U_2 \times ... \times U_i) \bigr) \cap \Delta =
U_1 \times U_2 \times \cdots \times U_i$.
The second assertion now follows from Lemma \ref{lemma - action on std hb}.
\qed\vskip 3mm

\begin{rk}\rm
\label{rk - why Delta}
The triviality of the ${\rm Stab}_{\Delta}(x_i)$-action on ${\rm HB}(x_i)$ is the main
reason why we work with $\Delta$ instead of $\Gamma$: we do not have the same property for
$\Gamma_i$ since $H_0$ usually acts non-trivially on ${\rm HB}(x_i)$.
\end{rk}

\begin{lemma}
\label{lemma - fd}
The subtree $F := \Gamma_0.R$ is a fundamental domain for the $\Delta$-action on $T$.
\end{lemma}

{\it Proof.~}
We have: $\Delta.F=(\Delta.\Gamma_0).R=\Gamma.R=T$.
Now assume that there are $\delta \!\in\! \Delta$ and $g,h \!\in\! \Gamma_0$ and
$i, \, j \geq 0$ such that $\delta.(g.x_i) = h.x_j =: x$. Then first of all
$i = j = \ell(x)$. Secondly,
we may and shall assume $i>0$ (otherwise $\delta=1$ by Lemma \ref{lemma - Levi}).
So $h{}^{-1}{\delta}g$ is in $\Gamma_i \, \subseteq \, {\Delta}H_0$.
Since $h{}^{-1}{\delta}h \!\in\! \Delta$, we conclude that $h{}^{-1}g$ is also in ${\Delta}H_0$.
However, $\Gamma_0\cap\Delta.H_0=H_0$ by Lemma \ref{lemma - Levi},
so that $h{}^{-1}g \!\in\! H_0$. Therefore: $g.x_i = h.x_i = h.x_j$.
\qed\vskip 3mm

This lemma shows in particular that $\Delta$ acts simply transitively on the set of all
vertices of level $0$.
Using again the coset representatives $\{\gamma_1:=1$; $\gamma_2$; $\cdots$ ; $\gamma_k\}$
for $\Gamma_0/H_0$,
we see that $F$ is the union of the $k$ rays $R_s := \gamma_s.R$ with $1 \leq s \leq k$.
We number the vertices of $R_s$ by $x_{i,s}:= \gamma_s.x_i$ for all $i \geq 0$ (so
$x_{0,s} = x_0$ for all $s$ and $x_{i,1} = x_i$ for all $i$).

\begin{claim}
\label{cl - action on hs and hb}
For any $1 \leq s \leq k$ and $i>0$, we have:
\begin{enumerate}
\item[(i)] The group $\langle U_{j,s}\mid j>i \rangle$ acts simply
transitively on ${\rm HS}(x_{i,s})$.
\item[(ii)] The stabilizer ${\rm Stab}_{\Delta}(x_{i,s})$ is equal to
$U_{1,s} \times \cdots \times U_{i,s}$ and acts trivially on ${\rm HB}(x_{i,s})$.
\end{enumerate}
\end{claim}

{\it Proof.~} The first statement follows directly from Fact \ref{fact - action on std hs},
and the second from
${\rm Stab}_{\Delta}(x_{i,s}) = \gamma_s{\rm Stab}_{\Delta}(x_i)\gamma_s^{-1}$ (here we use
that $\Delta$ is normal in $\Gamma$) and Lemma \ref{lemma - action on std hs and hb}.
\qed

\begin{lemma}
\label{lemma - action on hs and hb}
For any vertex $x$ of level $>0$, we have:
\begin{enumerate}
\item[(i)] the group ${\rm Stab}_\Delta(x)$ acts trivially on ${\rm HB}(x)$;
\item[(ii)] if $\delta,\delta'$ are two elements in $\Delta$ with $\delta.x = \delta'.x$,
then $\delta \!\mid_{{\rm HB}(x)}=\delta' \!\mid_{{\rm HB}(x)}$.
\end{enumerate}
\end{lemma}

{\it Proof.~} (ii) is a straightforward consequence of (i), and (i) follows from combining
Claim \ref{cl - action on hs and hb} and Lemma \ref{lemma - fd}.
\qed\vskip 3mm

Given $x \! \in \! VT$ with $\ell(x) > 0$, we can now label the elements in
${\rm HS}(x)$ by elements of $\Delta$.

\begin{definition}
\label{def - labels on hs}
\begin{enumerate}
\item[(a)] For $x = x_{i,s} \!\in\! F$ with $1 \leq s \leq k$ and $i >0$
and for any $y \!\in\! {\rm HS}(x)$, we define $\delta_{x,y}$ to be the unique element in
$\langle U_{j,s}\mid j>i \rangle$ satisfying  $\delta_{x,y}.x=y$
(see Claim \ref{cl - action on hs and hb}).
\item[(b)] For arbitrary $x \! \in \! VT$ with $\ell(x)>0$ and $y \!\in\! {\rm HS}(x)$, we
choose $\delta \!\in\! \Delta$ such that $\delta.x \!\in\! F$ and set
$\delta_{x,y} := \delta{}^{-1}\delta_{\delta(x),\delta(y)}\delta$.
\end{enumerate}
\end{definition}

The definition in point (b) above makes sense thanks to:

\begin{claim}
\label{cl - d well-defined}
The above introduced element $\delta_{x,y}$ is well-defined.
\end{claim}

{\it Proof.~}
Assume that also $\delta'.x \!\in\! F$ for some $\delta' \! \in \! \Delta$
and set: $z:=\delta'.x=x_{i,s} = \delta.x$ (see Lemma \ref{lemma - fd}).
Then $\delta'' := \delta'\delta{}^{-1} \!\in\! {\rm Stab}_\Delta(z) =
U_{1,s} \times ... \times U_{i,s}$ .
So firstly, $\delta''$ acts trivially on ${\rm HS}(z) =
{\rm HS}(\delta.x) = \delta.{\rm HS}(x)$ by Lemma \ref{lemma - action on hs and hb},
hence  $\delta'.y=\delta''.(\delta.y)=\delta.y$; and secondly $\delta''$ commutes with
$\langle U_{j,s}\mid j>i \rangle$,  in particular with $\delta_{z,\delta(y)}$.
(Here we use again that $\Gamma$ is of {\it directly} split type.)
Therefore
${\delta'}{}^{-1}\delta_{\delta'(x),\delta'(y)}\delta'
= \delta{}^{-1}\delta''{}^{-1}\delta_{z,\delta(y)}\delta''\delta
= \delta{}^{-1}\delta_{\delta(x),\delta(y)}\delta$.
\qed\vskip 3mm

Having established this, it is easy to check the following calculation rules:

\begin{lemma}
\label{lemma - d's}
For $x \! \in \! VT$ with $\ell(x)>0$ and $y \!\in\! {\rm HS}(x)$, we have:
\begin{enumerate}
\item[(i)] $\delta_{x,y}$ is in $\Delta$ and $\delta_{x,y}.x=y$;
\item[(ii)] $h\delta_{x,y}h{}^{-1}=\delta_{h(x),h(y)}$ for all $h \!\in\! \Delta$;
\item[(iii)] $\delta_{y,x}=\delta_{x,y}{}^{-1}$;
\item[(iv)] $\delta_{y,z}\delta_{x,y} = \delta_{x,z}$ for all $z \! \in \! {\rm HS}(x)$.
\end{enumerate}
\end{lemma}

{\it Proof.~}
(i) is clear by definition.
To prove (ii), we note that if $\delta \!\in\! \Delta$ is such that
$\delta.x \!\in\! F$, then we can choose ${\delta}h^{-1} \!\in\! \Delta$
in the definition of $\delta_{h(x),h(y)}$ in order to have
$({\delta}h^{-1}).(h.x)=\delta.x \!\in\! F$.
Therefore, applying part (b) of Definition \ref{def - labels on hs}, we obtain:
$\delta_{h(x),h(y)} =
({\delta}h^{-1})^{-1}\delta_{\delta(x),\delta(y)}({\delta}h^{-1}) =
h\delta_{x,y}h{}^{-1}$.

\vskip 3mm

According to (ii), it is enough to prove (iii) and (iv) for $x \!\in\! F$.
So assume $x = x_{i,s}$ with $i > 0$ and $1 \leq s \leq k$. Set
$\delta:=\delta_{x,y}{}^{-1}$, and observe that this is an element of
$\langle U_{j,s}\mid j>i \rangle$ sending $y$ to  $x \!\in\! F$ as well as
$\delta.y = x$ to $\delta.x$. So by part (a) of Definition \ref{def - labels on hs},
$\delta = \delta_{\delta(y),\delta(x)}$.
Therefore, now by part (b) of  Definition \ref{def - labels on hs},
$\delta_{y,x} = \delta^{-1}\delta_{\delta(y),\delta(x)}\delta =
\delta{}^{-1}\delta\delta = \delta = \delta_{x,y}{}^{-1}$, proving (iii).
In order to verify (iv), we first note that, again by part (b) of  Definition
\ref{def - labels on hs},
$\delta_{y,z} = \delta{}^{-1}\delta_{\delta(y),\delta(z)}\delta =
\delta{}^{-1}\delta_{x,\delta(z)}\delta$, which is an element of
$\langle U_{j,s}\mid j>i \rangle$. Hence also  $\delta_{y,z}\delta_{x,y}$ is an
element of $\langle U_{j,s}\mid j>i \rangle$, and by (i), it sends $x$ to $z$.
Therefore, by part (a) of  Definition \ref{def - labels on hs},
$\delta_{y,z}\delta_{x,y} = \delta_{x,z}$, proving (iv).
\qed\vskip 3mm

\subsection{Connected components of bounded level}
\label{ss - connected components}
We are now going to generalize the strategy of \cite[Section 8.3]{BurMoz} to the
present situation. We recall the definition of a certain graph $\hbox{\cal G}_i$ associated
to $T$ and a given level $i > 0$.

\begin{definition}
\label{def - graph}
For $i > 0$, we denote by $T_i$ the subforest of $T$ spanned by all vertices of level
$\leq i$. For any vertex $x$ with $\ell(x) \leq i$, we denote by $C_i(x)$ the
connected component of $T_i$ containing $x$. We set
$${\cal C}_i := \{ C_i(x) \mid x \! \in \! VT \ {\rm and} \ \ell(x) \leq i \},$$
and define a (simplicial) graph
$$\hbox{\cal G}_i = (V\hbox{\cal G}_i,E\hbox{\cal G}_i)
\quad {\rm with} \ {\rm vertex} \ {\rm set} \quad V\hbox{\cal G}_i = {\cal C}_i$$
by declaring a two element subset $\{ X, Y \}$ of  $V\hbox{\cal G}_i$ to be in the edge set
$E\hbox{\cal G}_i$ if and only if there exist vertices $x \! \in \! X$, $y \! \in \! Y$ such
that $\ell(x) = \ell(y) = i$ and $y \! \in \! {\rm HS}(x)$.
\end{definition}

\begin{rks} \rm
\label{rks - graph}
1. If  $\{ X, Y \} \! \in \! E\hbox{\cal G}_i$, the vertices $x \! \in \! X$ and
$y \! \in \! Y$ in the above definition are uniquely determined:
they are the unique terminal vertices of $X$ and $Y$, respectively, contained in any
geodesic in $T$ which starts in a vertex of $X$ and ends in a vertex of $Y$.

2. If $X,Y \! \in \! {\cal C}_i$ are given, there is a unique geodesic in
$\hbox{\cal G}_i$ connecting them. It is obtained from an arbitrary geodesic $p$ in $T$ which
connects a vertex of $X$ to a vertex of $Y$ by first intersecting $p$ with $T_i$,
yielding a disjoint union $p'$ of geodesics in $T_i$, then replacing these subgeodesics of
$p'$ with their corresponding connected components in ${\cal C}_i$ (considered as elements
of $V\hbox{\cal G}_i$) and connecting the adjacent ones among them by edges in $E\hbox{\cal
G}_i$.

3. A path ${\cal P}$ in $\hbox{\cal G}_i$, i.e. a finite sequence ${\cal
P}=(X_0,\ldots,X_n)$ with $X_j \! \in \! {\cal C}_i$ for all $j \leq n$ and $\{ X_j, X_{j+1}
\} \! \in \! E\hbox{\cal G}_i$ for all $j < n$, is a geodesic in $\hbox{\cal G}_i$ if and
only if $X_{j+2}$ is neither equal nor adjacent to $X_j$ for all $j \leq n - 2$. In fact,
there is only one path ${\cal P}$ between $X_0$ and $X_n$ in $\hbox{\cal G}_i$ with this last
mentioned property. This follows from the uniqueness of geodesics in $T$.
\end{rks}

The elements $\delta_{x,y}$ introduced in the previous Subsection now provide us with
transformations $\tau_{X,Y}$ mapping $X \! \in \! {\cal C}_i$ onto $Y \! \in \! {\cal C}_i$.

\begin{definition}
\label{def - tau}
Let $X,Y \! \in \! {\cal C}_i$ be given.
\begin{enumerate}
\item[(a)] If $X$ and $Y$ are adjacent in $\hbox{\cal G}_i$, and if $x \! \in VX$,
$y \! \in VY$ are such that $y \! \in \! {\rm HS}(x)$, we set
$\tau_{X,Y}:=\delta_{x,y}$. This is well-defined by Remark \ref{rks - graph} (1).
\item[(b)] If $X$ and $Y$ are arbitrary, we consider the unique geodesic
$(X=X_0$, $X_1$ , $\cdots$, $X_n=Y)$ from $X$ to $Y$ in $\hbox{\cal G}_i$, and set
$\tau_{X,Y} := \tau_{X_{n-1},X_n}\tau_{X_{n-2},X_{n-1}}\cdots\tau_{X_0,X_1}$. We also set
$\tau_{X,X} := 1$.
\end{enumerate}
\end{definition}

\begin{rk}\rm
\label{rk - path} We note that
$\tau_{X,Y} = \tau_{Z_{m-1},Z_m}\tau_{Z_{m-2},Z_{m-1}}\cdots\tau_{Z_0,Z_1}$ is
true for {\it any} path ${\cal P}=(X=Z_0$, $Z_1$ , $\cdots$, $Z_m=Y)$ connecting $X$ and
$Y$ in $\hbox{\cal G}_i$. Indeed, if  ${\cal P}$ is not a geodesic, then by Remark
\ref{rks - graph} (3), there exists a subpath $(A,B,C)$ of ${\cal P}$ such that $A = C$
or $\{ A, C \} \! \in \! E\hbox{\cal G}_i$. In the first case, we delete $B$ and $C = A$ from
${\cal P}$ and observe that $\tau_{B,A}\tau_{A,B} = 1$ by Lemma \ref{lemma - d's} (iii).
In the second case, we delete $B$ from ${\cal P}$ and note that
$\tau_{B,C}\tau_{A,B} = \tau_{A,C}$ by Lemma \ref{lemma - d's} (iv).
So simultaneously shortening the path ${\cal P}$ and the product of $\tau$'s associated
to it without changing the value of this product, we finally obtain the geodesic
between $X$ and $Y$ in $\hbox{\cal G}_i$, and the associated product is precisely the one
occurring in the above definition of $\tau_{X,Y}$.
\end{rk}

\begin{lemma}
\label{lemma - tau's}
For any $X, Y, Z \!\in\! {\cal C}_i$, we have:
\begin{enumerate}
\item[(i)] $\tau_{X,Y} \!\in\! \Delta$ and $\tau_{X,Y}(X)=Y$;
\item[(ii)] $h\tau_{X,Y}h{}^{-1} = \tau_{h(X),h(Y)}$ for all $h \!\in\! \Delta$;
\item[(iii)] $\tau_{Y,X} = \tau_{X,Y}{}^{-1}$;
\item[(iv)] $\tau_{Y,Z}\tau_{X,Y} = \tau_{X,Z}$.
\end{enumerate}
\end{lemma}

{\it Proof.~}
The properties (i) to (iii) follow directly from the corresponding statements in
Lemma \ref{lemma - d's}, and (iv) is a consequence of Remark \ref{rk - path}.
\qed\vskip 3mm

We close this subsection by introducing truncated versions of $T$, $R$, $F$ and $\Delta$.
For $i > 0$, we set $Y_i := C_i(x_0)$, $R^{(i)} := R \cap Y_i$,
$F^{(i)} := F \cap Y_i = \Gamma_0.R^{(i)}$ and
$$\Delta_i := \langle U_{j,s} \mid 1 \leq j \leq i \quad {\rm and} \quad 1 \leq s \leq k
\rangle .$$

\begin{lemma}
\label{lemma - uniform connected component}
\begin{enumerate}
\item[(i)] $R^{(i)}$ is a fundamental domain for the action of
$\langle \Gamma_0 \cup \Gamma_i \rangle = {\rm Stab}_{\Gamma}(Y_i)$ on $Y_i$;
\item[(ii)] $F^{(i)}$ is a fundamental domain for the action of
$\Delta_i = {\rm Stab}_{\Delta}(Y_i)$ on $Y_i$;
\item[(iii)]$\Delta_i$ (or more precisely: its image in ${\rm Aut}(Y_i)$) is a
uniform $Y_i$-lattice.
\end{enumerate}
\end{lemma}

{\it Proof.~} Set $\Gamma' := \langle \Gamma_0 \cup \Gamma_i \rangle$, which is obviously
contained in ${\rm Stab}_{\Gamma}(Y_i)$. It is clear that $\Gamma'.R^{(i)}$ contains a
neighborhood of $R^{(i)}$ in $Y_i$. Thus $\Gamma'.R^{(i)}$ is open in $Y_i$. But as a
subgraph of $Y_i$, it is also closed. Since $Y_i$ is connected, we conclude
$\Gamma'.R^{(i)} = Y_i$. (Side remark: this kind of reasoning already occurs in
\cite{Ser}.) Thus $\Gamma' \setminus Y_i = R^{(i)}$. If $\gamma \! \in \! \Gamma$
stabilizes $Y_i$, then $\gamma.x_j \! \in \! \Gamma'.x_j$ for any $j \leq i$, and since
${\rm Stab}_{\Gamma}(x_j) = \Gamma_j \subseteq \Gamma'$, we obtain
$\gamma \! \in \! \Gamma'$. This proves (i).

\vskip 3mm

The proof of (ii) is similar. We just observe that also $\Delta_i.F^{(i)}$ (which
contains a neighborhood of $x_0$ by definition of $F^{(i)}$) is open in $Y_i$ due to the
transitivity properties of the $U_{j,s}$ (see Fact \ref{fact - transitivity of root gp}).
This implies $\Delta_i.F^{(i)} = Y_i$ and hence $\Delta_i \setminus Y_i = F^{(i)}$.
Since ${\rm Stab}_{\Delta}(x) \subseteq \Delta_i$ for all vertices $x \! \in \! F^{(i)}$
(Claim \ref{cl - action on hs and hb} (ii)), we also obtain $\Delta_i =
{\rm Stab}_{\Delta}(Y_i)$. Finally, (iii) is an immediate consequence of (ii).
\qed

\subsection{Permuting elements of the same level}
\label{ss - permuting elements of the same level}
We now replace the analysis of \og local data\fg which were used in \cite{BurMoz} and
\cite{Moz} by a more group theoretic construction. The latter is based on the observation
Lemma \ref{lemma - action on hs and hb} (ii) above and similar to the definition
of the $\delta_{x,y}$. However, one has to be a bit careful here since $\Delta$ does
not act transitively on the set of vertices of a given level.

\begin{definition}
For any given vertex $x$ of level $i>0$, let $s = s(x)$ be the unique integer with
$1 \leq s \leq k$ and $x \!\in\! \Delta.R_s$ (see Lemma \ref{lemma - fd}).
\begin{enumerate}
\item[(a)]
We choose once and for all a $\delta_x \!\in\! \Delta$ such that
$x=\delta_x.x_{i,s}=(\delta_x\gamma_s).x_i$, and we set:
$\gamma_x:=\delta_x\gamma_s$.
\item[(b)] For any two vertices $x,y$ of the same level $i>0$, we define:
$\gamma_{x,y}:=\gamma_y\gamma_x{}^{-1}$.
\end{enumerate}
\end{definition}

Note that we work with the fixed system of coset representatives
$\{\gamma_1:=1$; $\gamma_2$; $\cdots$ ; $\gamma_k\}$ for $\Gamma_0/H_0$ and that
$\gamma_{x,y}$ need not coincide with $\delta_{x,y}$ in case $x$ and $y$ are in the
same horosphere (see, however, (iv) below which we will need to prove
Claim \ref{claim - D in L_i}).

\begin{lemma}
\label{lemma - g's}
For vertices $x,y,z$ of the same level $i>0$, we have:
\begin{enumerate}
\item[(i)] $\gamma_{x,y}$ is in $\Gamma$ and $\gamma_{x,y}.x = y$;
\item[(ii)] $\gamma_{y,x}=\gamma_{x,y}{}^{-1}$;
\item[(iii)] $\gamma_{y,z}\gamma_{x,y}= \gamma_{x,z}$;
\item[(iv)] if $y \!\in\! \Delta.x$ then $\gamma_{x,y} \!\in\! \Delta$.
\end{enumerate}
\end{lemma}

{\it Proof.~}
(i), (ii) and (iii) are obvious by definition. If $x \!\in\! \Delta.R_s$ and
$y \!\in\! \Delta.R_t$ with $1 \leq s, \, t \leq k$, then
$\gamma_{x,y} \!\in\! \Delta\gamma_t\gamma_s{}^{-1}\Delta = \Delta\gamma_t\gamma_s{}^{-1}$
because $\Delta \triangleleft \Gamma$ (Lemma \ref{lemma - Levi}). This implies (iv).
\qed\vskip 3mm

The next property is a bit more subtle. It will be needed in order to obtain a reasonable
definition of the groups $L_i$ (see Definition \ref{def - Li} and Remark \ref{rks - Li} (1))
which are essential for the extensions to be constructed in the next subsection.

\begin{lemma}
\label{lemma - restriction of g's}
If $x,y$ are vertices of level $i > 0$,
$x' \!\in\! {\rm HS}(x)$ and $y':= \gamma_{x,y}(x')$, then we have:
$\gamma_{x,y}\!\mid_{{\rm HB}(x)}=\gamma_{x',y'}\!\mid_{{\rm HB}(x)}$.
\end {lemma}

{\it Proof.~}
By Lemma \ref{lemma - g's} (iii), we have:
$\gamma_{x',y'}=\gamma_{y,y'}\gamma_{x,y}\gamma_{x',x}$.
Since $x' \!\in\! {\rm HS}(x)$ and hence
$y' \!\in\! \gamma_{x,y}.{\rm HS}(x)={\rm HS}(y)$, $x' \! \in \Delta.x$ and
$y' \! \in \! \Delta.y$ by the transitive action of $\Delta$ on ${\rm HS}(x)$,
respectively ${\rm HS}(y)$ (combine Claim \ref{cl - action on hs and hb}
and Lemma \ref{lemma - fd}). Thus Lemma \ref{lemma - g's} (iv) tells us that
$\gamma_{x',x}$ and $\gamma_{y,y'}$ are in $\Delta$. Therefore
$\gamma_{x',y'} \!\in\! \Delta\gamma_{x,y}\Delta = \gamma_{x,y}\Delta$, the last
equality again following from $\Delta \triangleleft \Gamma$. Hence there exists a
$\delta \!\in\! \Delta$ such that $\gamma_{x',y'}=\gamma_{x,y}\delta$. Since
$\delta.x' = \gamma_{x,y}{}^{-1}\gamma_{x',y'}.x'= \gamma_{x,y}{}^{-1}(y') =  x'$, we
have $\delta \!\in\! {\rm Stab}_\Delta(x')$. Therefore $\delta$ acts trivially on
${\rm HB}(x') = {\rm HB}(x)$ (Lemma \ref{lemma - action on hs and hb}), and finally:
$\gamma_{x',y'}\!\mid_{{\rm HB}(x)}=\gamma_{x,y}\delta\!\mid_{{\rm HB}(x)}=
\gamma_{x,y}\mid_{{\rm HB}(x)}$.
\qed\vskip 3mm

\subsection{Extensions and commensurators}
\label{ss - extensions and commensurators}
We can now define suitable extensions of isomorphisms between subtrees $X,Y \! \in \!
{\cal C}_i$ of $T$ (see Proposition \ref{prop - unique extension} below). For this we need
to define certain subgroups of $L$.

\begin{definition}
\label{def - Li}
Let $i>0$ be a positive integer.
By $L_i$ we denote the group of all level-preserving automorphisms $h \!\in\! {\rm Aut}(T)$
satisfying:
\begin{enumerate}
\item[(a)] for any vertex $x$ with $\ell(x) = i$,
$h\!\mid_{{\rm HB}(x)} = \gamma_{x,h(x)} \!\mid_{{\rm HB}(x)}$;
\item[(b)]  for all $X,Y \! \in \! {\cal C}_i$,
$(h\tau_{X,Y}h{}^{-1})\!\mid_{h(X)}=\tau_{h(X),h(Y)} \!\mid_{h(X)}$.
\end{enumerate}
\end{definition}

\begin{rks}\rm
\label{rks - Li}
1. $L_i$ is indeed a subgroup of $L$. If $g,h \! \in \! L_i$ then $gh,h^{-1}$
clearly satisfy (b), and it follows from Lemma \ref{lemma - g's} that they also satisfy
(a). In the extension procedure, we only need property (a) of the above definition at certain
vertices $x$. However, one would not get a reasonable definition of $L_i$ if one did not
require (a) for all level $i$ vertices $x$.
But then one has to check that the different vertices in ${\rm HS}(x)$ do not lead to
different conditions for $h\!\mid_{{\rm HB}(x)}$: this is precisely what we established in
Lemma \ref{lemma - restriction of g's}.

2. It is sufficient to require Condition (b) for a {\it fixed} $X_0 \! \in \! {\cal C}_i$
(and all $Y \! \in \! {\cal C}_i$). This follows from Lemma \ref{lemma - tau's}:
$\tau_{X,Y} = \tau_{X_0,Y}\tau_{X_0,X}{}^{-1}$ for all $X,Y \! \in \! {\cal C}_i$.
\end{rks}

\begin{claim}
\label{claim - D in L_i}
For all $i>0$, the group $\Delta$ is contained in $L_i$.
\end{claim}

{\it Proof.~}
Condition (a) for elements of $\Delta$ follows from combining
Lemma \ref{lemma - action on hs and hb} (ii)  and Lemma \ref{lemma - g's} (iv).
Condition (b) follows from Lemma \ref{lemma - tau's} (ii).
\qed\vskip 3mm

Up to replacing local data by the elements $\gamma_{x,y}$, the next proposition is
similar to \cite[Proposition 8.5]{BurMoz}.

\begin{prop}
\label{prop - unique extension}
For any two $X,Y \! \in \! {\cal C}_i$ and any level-preserving isomorphism
$h: X \to Y$, there is a unique extension of $h$ to an automorphism $E(h) \!\in\! L_i$.
\end{prop}

{\it Proof.~}
We first show how we can (and must) extend $h$ to horoballs and elements of ${\cal C}_i$
which are neighboring $X$.

\vskip 3mm

(1) {\it Definition of $E(h)$ on ${\rm HB}(x)$.}
According to Condition (a) of Definition \ref{def - Li}, we have only one choice to define
$E(h)$ on ${\rm HB}(x)$ for $x \! \in \! VX$ with $\ell(x) = i$, namely by
$E(h)\!\mid_{{\rm HB}(x)} = \gamma_{x,y}\!\mid_{{\rm HB}(x)}$, where $y := h(x)$.

\vskip 3mm

(2) {\it Definition of $E(h)$ on $Z$.}
For $x \neq z \! \in \! {\rm HS}(x)$, we set $Z := C_i(z) \! \in \! {\cal C}_i$,
$z' := E(h)(z) = \gamma_{x,y}(z)$ and $Z' := C_i(z')$. Note that $z' \! \in \! {\rm HS}(y)$.
According to Condition (b) of Definition \ref{def - Li}, we have to define $E(h)$ on $Z$ by
$E(h)\!\mid_Z = \tau_{Y,Z'}h\tau_{Z,X}\!\mid_Z$.

\vskip 3mm

(3) {\it Equality of the two definitions on $\{ z \} = Z \cap {\rm HB}(x)$.}
Since $z' \! \in \! {\rm HS}(y)$ and $z \! \in \! {\rm HS}(x)$, Definition \ref{def - tau}
implies $\tau_{Y,Z'}(y) = z'$ and $\tau_{Z,X}(z) = x$. Hence
$\tau_{Y,Z'}h\tau_{Z,X}\!\mid_Z(z)=\tau_{Y,Z'}h(x)=\tau_{Y,Z'}(y)=z'= \gamma_{x,y}(z)$.

\vskip 3mm

(4) {\it Global definition of $E(h)$.}
Now an easy induction along geodesics in $\hbox{\cal G}_i$ which start in $X$ shows that we can, in
precisely one way, extend the definition of $E(h)$ to horoballs and elements of ${\cal C}_i$
in the neighborhood of ${\tilde Z}$ for any ${\tilde Z} \! \in \! {\cal C}_i$ on which $E(h)$
is already defined. We thus get an automorphism $E(h)$ of $T$, and it is clear by construction
that $E(h)$ is level-preserving.

\vskip 3mm

(5) {\it $E(h)$ is an element of $L_i$.}
For any horoball $H = {\rm HB}(v)$, ($v \! \in \! VT$ with $\ell(v) = i$), there is by
construction of $E(h)$ a $z_1 \! \in \! HS = {\rm HS}(v)$ such that
$E(h)\!\mid_H = \gamma_{z_1,z_2}\!\mid_H$, where $z_2 = E(h)(z_1)$. But then by Lemma
\ref{lemma - restriction of g's}, $E(h)\!\mid_H = \gamma_{w_1,w_2}\!\mid_H$ with
$w_2 = E(h)(w_1)$ for all $w_1 \! \in \! HS$, yielding Condition (a) of Definition
\ref{def - Li}. To check Condition (b), we first consider an edge $\{ A,B \}$ in
$\hbox{\cal G}_i$ such that $A$ is nearer to $X$ than $B$, and set $A' = E(h)(A)$,
$B' = E(h)(B)$. By the inductive definition of $E(h)$, we have
$E(h)\tau_{A,B}E(h)^{-1}\!\mid_{A'} = \tau_{A',B'}\!\mid_{A'}$. Similarly for all edges
on the geodesic in $\hbox{\cal G}_i$ between $X$ and $B$. Applying Definition \ref{def - tau}
to $\tau_{X,B}$, we now obtain $E(h)\tau_{X,B}E(h)^{-1}\!\mid_{Y} = \tau_{Y,B'}\!\mid_{Y}$.
Since $B \! \in \! {\cal C}_i$ can be chosen arbitrarily, Remark \ref{rks - Li} (2) implies 
that Condition (b) is also satisfied.
\qed\vskip 3mm

The next proposition shows that the above extension procedure is compatible with commensurators.
This is \cite[Proposition 8.6]{BurMoz}, which we reproduce to correct minor misprints and
because we deal with the group $\Delta$ rather than $\Gamma$.
We denote the fact that a subgroup $N$ is of finite index in a group $H$ by writing:
$N <_{\rm f.i.}H$.
Since $\Delta <_{\rm f.i.} \Gamma$, the commensurators $C(\Delta)$ and $C(\Gamma)$
of these groups in $G$ coincide.
We keep the notations $Y_i = C_i(x_0)$ and $\Delta_i = {\rm Stab}_{\Delta}(Y_i)$ introduced 
in Lemma \ref{lemma - uniform connected component} and set $G_i := {\rm Aut}(Y_i)$,
$G_i' := \{ h \! \in \! G_i \mid \ell(h.x) = \ell(x)$ for all $x \! \in \! VY_i \}$.

\begin{prop}
\label{prop - extension-comm}
Let ${\overline \Delta_i}$ be the natural image of $\Delta_i$ in $G_i'$. Then for any
$g \! \in \! {\rm Comm}_{G_i'}({\overline \Delta_i})$, the unique extension 
$E(g) \! \in \! L_i$ provided by Proposition \ref{prop - unique extension} lies in the
commensurator ${\rm Comm}_L(\Delta) = C(\Delta) \cap L$. 

\end{prop}

{\it Proof.~}
We first remark that the map $E: G_i' \rightarrow L_i, \ h \mapsto E(h)$ is a 
homomorphism by the uniqueness statement of Proposition \ref{prop - unique extension}.
It is enough to show that for any $g \!\in\! {\rm Comm}_{G_i'}({\overline \Delta_i})$, 
we have $\Delta \cap E(g)\Delta E(g)^{-1} <_{\rm f.i.}\Delta$.
(Indeed, apply it to $g^{-1}$ and conjugate by $E(g)$ to deduce that
$\Delta \cap E(g)\Delta E(g)^{-1} <_{\rm f.i.} E(g) \Delta E(g)^{-1}$.)
By definition of a commensurator, there is
$\{\delta_j \}_{j=1}^r \subset \Delta_i$ such that, with
${\overline \delta_j} := \delta_j \! \mid_{Y_i}$,
${\overline \Delta_i}=
\bigsqcup_{j=1}^r {\overline \delta_j}({\overline \Delta_i}\cap g{\overline
\Delta_i}g^{-1})$. It is enough to prove that we have:
$\Delta=\bigcup_{j=1}^r \delta_j(\Delta \cap E(g)\Delta E(g)^{-1})$.
Let $\delta \!\in\! \Delta$.
Since $\tau_{\delta.Y_i,Y_i} \circ \delta$ stabilizes $Y_i$, there is an index
$j$ such that
${\overline \delta_j}^{-1} \circ(\tau_{\delta.Y_i, Y_i} \circ \delta)\!\mid_{Y_i}
\in {\overline \Delta_i}\cap g{\overline \Delta_i}g^{-1}$.
In particular,
$\sigma :=  g^{-1} \circ \bigl( \delta_j^{-1} \circ \tau_{\delta.Y_i, Y_i} \circ \delta
\bigr) \! \mid_{Y_i} \circ g $ lies in  ${\overline \Delta_i}$.
Now we consider $E(g)^{-1}(\delta_j^{-1}\delta)E(g) \! \in \! L_i$ (recall that
$\Delta \subseteq L_i$ by Claim \ref{claim - D in L_i}) and prove
that its restriction to $Y_i$ coincides with the restriction of an element of $\Delta$
to $Y_i$. By uniqueness (Proposition \ref{prop - unique extension}), this will show that
$E(g)^{-1}(\delta_j^{-1}\delta)E(g) \!\in\! \Delta$.
The element $\delta$ being arbitrary in $\Delta$, this will finally prove:
$\Delta \cap E(g)\Delta E(g)^{-1}<_{\rm f.i.}\Delta$.
We are thus led to computing:

\vskip 3mm

$\bigl( E(g)^{-1}(\delta_j^{-1}\delta)E(g) \bigr)\!\mid_{Y_i}
=E(g)^{-1} \! \mid_{\delta_j^{-1}\delta.Y_i} \circ \, \delta_j^{-1}
\! \mid_{\delta.Y_i} \circ \, \delta \! \mid_{Y_i} \circ \, g
=E({\overline \delta_j} g)^{-1} \! \mid_{\delta.Y_i} \circ \, \delta \! \mid_{Y_i} \circ
\, g$.

\vskip 3mm

We make $\sigma$ appear on the right to obtain an element of ${\overline \Delta_i}$:

\vskip 3mm

\hskip 10mm $= E({\overline \delta_j} g)^{-1} \! \mid_{\delta.Y_i} \circ \,
\bigl(
(\tau_{\delta.Y_i,Y_i}\! \mid_{\delta.Y_i})^{-1}  \circ {\overline \delta_j} \circ \, g
\bigr)
\circ
\bigl( g^{-1} \circ {\overline \delta_j}^{-1} \circ \, \tau_{\delta.Y_i,Y_i}\!
\mid_{\delta.Y_i} \bigr) \circ
\delta \! \mid_{Y_i} \circ \, g$

\hskip 10mm $= E({\overline \delta_j} g)^{-1} \! \mid_{\delta.Y_i} \circ \,
\bigl(
(\tau_{\delta.Y_i,Y_i}\! \mid_{\delta.Y_i})^{-1}  \circ {\overline \delta_j} \circ \, g
\bigr)
\circ \, \sigma$

\hskip 10mm $=
\bigl( E({\overline \delta_j} g)^{-1} \! \mid_{\delta.Y_i} \circ \,
\tau_{Y_i,\delta.Y_i}\! \mid_{Y_i} \circ \, E({\overline \delta_j} g)\! \mid_{Y_i} 
\bigr) \circ \, \sigma$.

Since $E({\overline \delta_j}g) \! \in \! L_i$, this last product
is $\tau_{Y_i,Z} \! \mid_{Y_i} \! \circ \, \sigma$, where
$Z:=E({\overline \delta_j}g)^{-1}(\delta.Y_i)$. Because $\tau_{Y_i,Z} \! \in \! \Delta$
and $\sigma \! \in \! {\overline \Delta_i}$, the restriction of 
$E(g)^{-1}(\delta_j^{-1}\delta)E(g)$ to $Y_i$ coincides with the restriction of an
element of $\Delta$ to $Y_i$.
\qed\vskip 3mm

\subsection{Density}
\label{ss - density}
We are now in a position to complete the proof of our main density theorem. The main
tools will be Proposition \ref{prop - extension-comm} and Liu's theorem on the
arithmeticity of uniform tree lattices.

\vskip 3mm

{\it Proof of Theorem \ref{th - density}.}
Given any level-preserving isomorphism $\phi: T_1 \rightarrow T_2$ between two
finite subtrees $T_1$ and $T_2$ of $T$, we have to find an element
$h \! \in \! C(\Gamma) \cap L = C(\Delta) \cap L$ such that $h{\mid}_{T_1} = \phi$.
We choose $i > 0$ big enough so that the subtree $Y_i = C_i(x_0)$ of $T$ introduced in
Lemma \ref{lemma -  uniform connected component} satisfies the following two conditions:

\begin{enumerate}
\item[(i)]  the subtree $Y_i$ contains $T_1 \cup T_2$;
\item[(ii)] the tree $Y_i$ is not biregular.
\end{enumerate}

We briefly explain why Condition (ii) can always be achieved and why we need it.
First recall that the vertices of level $j < i$ in $Y_i$ have the same degree in $Y_i$ as
in $T$, namely $q_j + 1$ (Facts \ref{fact - valency at level 0} and \ref{fact - valency at
level > 0}). And the vertices of level $i$ are of degree $q_i$ in $Y_i$. So either $T$ is
not biregular, and there exists an $l$ such that $q_l \neq q_{l+2}$; then we choose
$i > l + 2$. Or $T$ is biregular, with degrees $q_0 + 1$ and $q_1 + 1$. In this case
$Y_i$ is not biregular for any $i \geq 2$. If $Y_i$ is not biregular, then we know from
Lemma \ref{lemma - L equals G}, Remark \ref{rk - biregular} and Lemma
\ref{lemma -  uniform connected component} (i) that all automorphisms of $Y_i$ are
automatically level--preserving. This means $G_i' = G_i = {\rm Aut}(Y_i)$ in the notation
of Proposition \ref{prop - extension-comm}.

\vskip 3mm

Now  ${\overline \Delta_i}$, the image of $\Delta_i$ in $G_i$, is a uniform
$Y_i$-lattice (Lemma \ref{lemma -  uniform connected component} (iii)).
So by Liu's theorem, the commensurator ${\rm Comm}_{G_i}({\overline \Delta_i})$ is dense in
$G_i = G_i'$. Hence there exists a $g \! \in \! {\rm Comm}_{G_i'}({\overline \Delta_i})$ 
satisfying $g{\mid}_{T_1} = \phi$. By Proposition \ref{prop - unique extension},
$g$ can be (uniquely) extended to an element $h:=E(g)$ of $L_i$. Then Proposition
\ref{prop - extension-comm} implies that  $h \! \in \! C(\Delta) \cap L$, and by 
construction $h{\mid}_{T_1} = \phi$.
\qed\vskip 3mm

Theorem \ref{th - density} has some immediate consequences which are worth to be mentioned.
The first one is the generalization of the example discussed in \cite[Section 8.3]{BurMoz}
to the class of all lattices of directly split Nagao type with associated non-biregular
tree $T$. Recall that $L = G$ in these cases (see Lemma \ref{lemma - L equals G}).

\begin{cor}
If $\Gamma$ is a $T$-lattice of directly split Nagao type and $T$ is not
biregular, then $C(\Gamma)$ is dense in $G = {\rm Aut}(T)$.
\qed\end{cor}

The next corollary combines Theorem \ref{th - density} and Proposition
\ref{prop - maximality of L}. It was stated as Theorem \ref{th - biregular density}
in the introduction. 

\begin{cor}
If $\Gamma$ is a $T$-lattice of directly split Nagao type and $T$ is
biregular, then either $C(\Gamma)$ is a subgroup of $L$ which is then dense in $L$ or
else the closure of $C(\Gamma)$ in $G$ contains $G^\circ$.
\qed\end{cor}

For regular trees, this has the following immediate consequence.

\begin{cor}
If $\Gamma$ is a $T$-lattice of directly split Nagao type such that $T$
is regular and $C(\Gamma)$ contains a type-interchanging automorphism, then $C(\Gamma)$
is dense in $G$.
\qed\end{cor}

If $T$ is biregular, it is in general hard to see (and we do not know of any general
method how to decide this question) whether $C(\Gamma)$ is included in $L$ or not.
However, there is one class of lattices of Nagao type where
$C(\Gamma) \not\subseteq L$ is obvious,
and this class of lattices we are going to discuss in the next section.

\vskip 6mm

\section{Applications to Moufang twin trees}
\label{s - mtt}

We introduce Moufang twin trees and recall some decompositions available for
groups acting on these trees: this enables us to connect this theory to the previous
section.
We also provide examples, from the Nagao lattice to more exotic
twin trees, including generalized rank 2 Kac-Moody groups.
General references for twin trees are \cite{RonTit1}, \cite{RonTit2}
and \cite[\S 9]{cdf}.

\subsection{Twin trees without groups}
\label{ss - ttt}
The notion of codistance relating two trees can be introduced without reference to any
group action \cite[\S 1]{RonTit1}.
In this first subsection, we introduce as many notions as possible purely in combinatorial
terms.
We will show later how they coincide with previously defined objects in the Nagao context,
when we have enough automorphisms.

\begin{definition}
\label{definition - codistance}
Let $(T_\pm,\delta^*)$ be a triple where $T_+$ and $T_-$ are trees and where $\delta^*$ is
a function with values in ${\bf N}$  on pairs of vertices of opposite signs in
$T_+ \sqcup T_-$.
We say that $(T_\pm,\delta^*)$ is a {\rm twin tree~}or, equivalently, that $
\delta^*$ is a {\rm codistance~} between $T_+$ and $T_-$ if $\delta^*$ satisfies:

\vskip 3mm

{\rm (Codist)~}For any $x_+ \!\in\! VT_+$ and any $y_- \!\in\! VT_-$,
we have  $\delta^*(y_-,x_+) = \delta^*(x_+,y_-)$, and furthermore,
setting $m:=\delta^*(x_+,y_-)$: 
for each $y_-' \!\in\! VT_-$ adjacent to $y_-$, one obtains
$\delta^*(x_+,y_-')=m \pm 1$;  if $m>0$ there is a unique $y_-'$ such that 
$\delta^*(x_+,y_-')=m + 1$;
and we require the similar conditions with $T_+$ and $T_-$ interchanged.

\vskip 3mm

In this case, we say that two vertices of opposite signs are {\rm opposite~}if their
codistance is $0$.

\vskip 3mm

An {\rm automorphism~} of the twin tree  $(T_\pm,\delta^*)$ is a pair 
$(\alpha_+, \alpha_-)$ of automorphisms of $T_+$, respectively $T_-$, such that
$\delta^*(\alpha_+(x_+),\alpha_-(y_-)) = \delta^*(x_+,y_-)$ for any $x_+ \! \in \! VT_+$
and any $y_- \! \in \! VT_-$. We denote the group of all automorphisms of 
$(T_\pm,\delta^*)$ by $A$.
\end{definition}

\vskip 3mm

Two vertices in a tree have the same type if they are at even distance from one another.
Adding that two vertices of opposite sign have the same type if they are at even
codistance, we obtain an equivalence relation (with two classes) on the vertices of $T_+
\sqcup T_-$ \cite[p. 406]{RonTit1}.
Two thick trees (i.e. two trees where all vertices have degree at least 3)
related by a codistance are isomorphic and biregular because any two
vertices of the same type have the same degree \cite[Proposition 1]{RonTit1}.
We denote the types by $0$ and $1$ and write them as subscripts when needed.
Note that we do not require automorphisms of twin trees to be type-preserving.
We denote by $A^{\circ}$ the subgroup of $A$ of all type-preserving automorphisms
of $(T_\pm,\delta^*)$. Obviously, $[A : A^{\circ}] \leq 2$.

\vskip 3mm

Let $x_0^+$ and $x_1^+$ be two adjacent vertices in $T_+$.
We can choose two vertices $x^-_0$ and $x^-_1$ such that $x^+_i$ and $x^-_i$ are
opposite, $i \!\in\! \{ 0;1 \}$: these vertices will be called the {\it standard
positive~}or {\it negative~}vertices of type $0$  or $1$, accordingly.
A pair of geodesic lines $L_\pm \subset T_\pm$ is called a {\it twin apartment~}
if any vertex in $L:=L_+ \cup L_-$  has a unique opposite in $L$.
By \cite[Proposition 3.5]{RonTit1} there is a unique twin apartment containing the four
standard vertices, which we call the {\it standard twin apartment}.
It follows from the axioms of a codistance that the codistance between to vertices of
opposite sign in a twin apartment is the distance between any of the two points to the
unique opposite of the other one in the twin apartment.

\vskip 3mm

A straightforward consequence of (Codist) is a monotonicity property
distinguishing a subfamily of  geodesic rays and boundary
points in each tree \cite[3.1-3.4]{RonTit1}: for each pair of non-opposite vertices $x_+$
and $y_-$ there is a unique ray in $T_+$ (resp. in $T_-$) emanating from $x_+$ (resp.
$y_-$) along which the codistance from $y_-$ (resp. $x_+$) is increasing.
This defines a unique point $\xi(x_+)_{y_-}$ of the ideal boundary $\partial_\infty T_+$ and
a unique point $\eta(y_-)_{x_+}$ of the ideal boundary $\partial_\infty T_-$, so that the rays
are $[x_+;\xi(x_+)_{y_-})$ and $[y_-;\eta(y_-)_{x_+})$.
Using a standard pair of opposite vertices $x_i^\pm$, we see that the above monotonicity
arguments  naturally define in $T_\pm$ a union of geodesic rays, which we call the {\it
standard infinite star of type $i$~} and which we denote by ${\rm St}^\infty(x_i^\pm)$.
The number of rays from $x_i^\pm$ in ${\rm St}^\infty(x_i^\pm)$ is the valency of this
vertex.

\subsection{Root groups and Moufang condition}
\label{ss - Moufang}
The root system $\Phi$ attached to an arbitrary  Coxeter system is defined in
 \cite[Sect. 5]{Tit87} or in \cite[II.5]{Hum90}.
We are interested in the specific case where the Coxeter group is the
infinite dihedral group $D_\infty={\bf Z}/2 * {\bf Z}/2$, say with generators $s_0$ and $s_1$.
The associated Coxeter complex is the tiling of the real line ${\bf R}$ by the segments
$[n;n+1]$, $n \!\in\! {\bf Z}$.
Let $s_0$ (resp. $s_1$) be the reflection $x \mapsto -x$ (resp. $x \mapsto 2-x$).
The {\it roots~}of $\Phi$ are the half-lines defined by the integers, the positive
ones  being those containing $[0,1]$.
Each vertex has type 0 or 1, the boundary of a root is called its {\it vertex}, and the {\it
type~}of a root $a$ is the type of its vertex.
Two roots $a$ and $b$ are {\it prenilpotent~}with one another if
$a \subseteq b$ or $b \subseteq a$ \cite{Tit87}.
Let $a_0:=[0;+\infty)$ and $a_1:=(-\infty;1]$ be the roots whose intersection is the edge 
$E:=[0;1]$.
A useful viewpoint on twin apartments is to see them as two Tits cones of the same Weyl
group glued along their tip \cite[5.3.2]{RemAst}.
In our case, a twin apartment is the double cone generated by the above tiling of the real
line. We recover the geodesics $L_\pm$ as the affinizations of the double
cone.
A {\it twin root~} generated by a root $a$ of any sign is then the half-plane bounded by
the line passing through the origin and the vertex of $a$, and containing $a$.
We denote it by $a$, too. Two twin roots generated by prenilpotent roots $a$ and $b$
are also called prenilpotent.

\vskip 3mm

For a single building, the {\it Moufang condition~} is discussed in detail in \cite[\S
6]{RonLec}.
For twin trees it is introduced in \cite[pp. 475-476]{RonTit1} as follows.
Let $a$ be a twin root in the twin apartment $L$.
Let us denote by $U_a$ the group of automorphisms of $(T_\pm,\delta^*)$ fixing the half
twin apartment $a$ and every edge having a vertex in the interior of $a$.
By \cite[Proposition 4.1]{RonTit1}, the group $U_a$ acts freely on the set of twin apartments
containing $a$, but it may be trivial.
If the action is transitive we say that $U_a$ is a {\it root group}.

\begin{definition}
\label{definition - moufang}
We say that a twin tree $(T_\pm,\delta^*)$ satisfies the {\rm Moufang property~} or is a
{\rm Moufang twin tree~} if there is a twin apartment $L$ in which $U_a$ is a root group for
each twin root $a \subset L$.
We denote by $\Lambda$ the subgroup of $A^{\circ}$ generated by the root groups $U_a$.
\end{definition}

The definition seems to depend on the choice of a twin apartment, but if a twin tree is
Moufang, then $\Lambda$ is transitive on the set of twin apartments \cite[Proposition
4.5]{RonTit1}.
Note that so far the roles of the trees $T_+$ and $T_-$ are symmetric.
The Moufang property says that the valency at a vertex $v$ is 1 greater than the order of 
the root group attached to a twin root whose boundary contains $v$, so the trees
are locally finite if and only if the root groups are finite, in which case the pointwise
fixator of a twin apartment, which we denote by $H$, is finite too.

\vskip 3mm

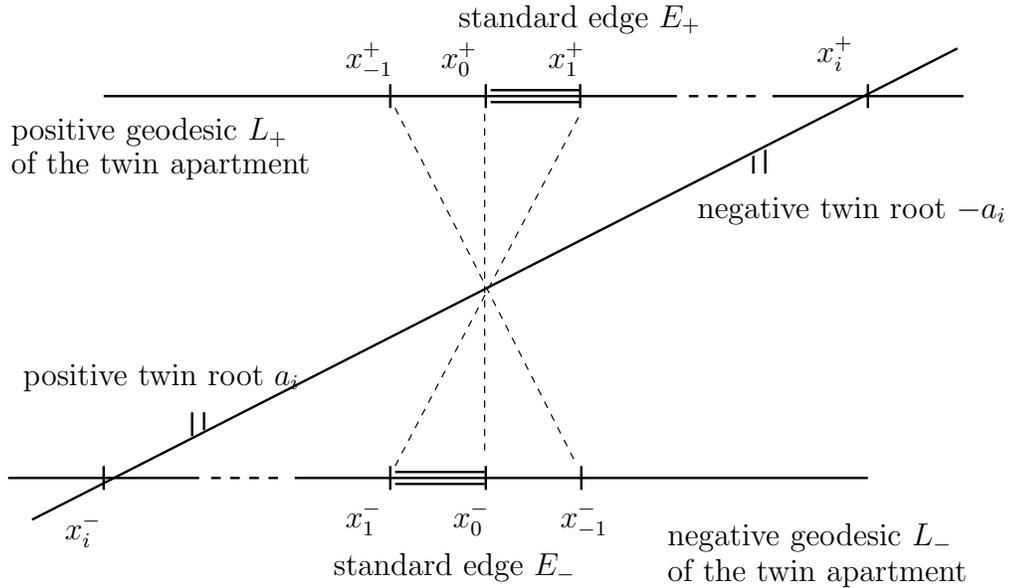
\begin{pict}
\label{pict - double Tits cone}
\end{pict}

\vskip 3mm
\centerline{\begin{picture}(0,0)%
\includegraphics{Picture3.pstex}%
\end{picture}%
\setlength{\unitlength}{1973sp}%
\begingroup\makeatletter\ifx\SetFigFont\undefined
\def\x#1#2#3#4#5#6#7\relax{\def\x{#1#2#3#4#5#6}}%
\expandafter\x\fmtname xxxxxx\relax \def\y{splain}%
\ifx\x\y   
\gdef\SetFigFont#1#2#3{%
  \ifnum #1<17\tiny\else \ifnum #1<20\small\else
  \ifnum #1<24\normalsize\else \ifnum #1<29\large\else
  \ifnum #1<34\Large\else \ifnum #1<41\LARGE\else
     \huge\fi\fi\fi\fi\fi\fi
  \csname #3\endcsname}%
\else
\gdef\SetFigFont#1#2#3{\begingroup
  \count@#1\relax \ifnum 25<\count@\count@25\fi
  \def\x{\endgroup\@setsize\SetFigFont{#2pt}}%
  \expandafter\x
    \csname \romannumeral\the\count@ pt\expandafter\endcsname
    \csname @\romannumeral\the\count@ pt\endcsname
  \csname #3\endcsname}%
\fi
\fi\endgroup
\begin{picture}(12066,7362)(57,-6535)
\put(301,-3961){\makebox(0,0)[lb]{\smash{\SetFigFont{12}{14.4}{rm}positive twin root $a_i$}}}
\put(151,-886){\makebox(0,0)[lb]{\smash{\SetFigFont{12}{14.4}{rm}positive geodesic $L_+$}}}
\put(151,-1306){\makebox(0,0)[lb]{\smash{\SetFigFont{12}{14.4}{rm}of the twin apartment }}}
\put(8401,-5986){\makebox(0,0)[lb]{\smash{\SetFigFont{12}{14.4}{rm}negative geodesic $L_-$}}}
\put(826,-5911){\makebox(0,0)[lb]{\smash{\SetFigFont{12}{14.4}{rm}$x_i^-$}}}
\put(10276, 89){\makebox(0,0)[lb]{\smash{\SetFigFont{12}{14.4}{rm}$x_i^+$}}}
\put(4351,-5761){\makebox(0,0)[lb]{\smash{\SetFigFont{12}{14.4}{rm}$x_1^-$}}}
\put(5701,-5761){\makebox(0,0)[lb]{\smash{\SetFigFont{12}{14.4}{rm}$x_0^-$}}}
\put(7051,-5761){\makebox(0,0)[lb]{\smash{\SetFigFont{12}{14.4}{rm}$x_{-1}^-$}}}
\put(5551, 14){\makebox(0,0)[lb]{\smash{\SetFigFont{12}{14.4}{rm}$x_0^+$}}}
\put(6901, 14){\makebox(0,0)[lb]{\smash{\SetFigFont{12}{14.4}{rm}$x_1^+$}}}
\put(4351, 14){\makebox(0,0)[lb]{\smash{\SetFigFont{12}{14.4}{rm}$x_{-1}^+$}}}
\put(8401,-6436){\makebox(0,0)[lb]{\smash{\SetFigFont{12}{14.4}{rm}of the twin apartment }}}
\put(4201,-6361){\makebox(0,0)[lb]{\smash{\SetFigFont{12}{14.4}{rm}standard edge $E_-$}}}
\put(5776,539){\makebox(0,0)[lb]{\smash{\SetFigFont{12}{14.4}{rm}standard edge $E_+$}}}
\put(8776,-1861){\makebox(0,0)[lb]{\smash{\SetFigFont{12}{14.4}{rm}negative twin root $-a_i$}}}
\end{picture}
}
\vskip 3mm

We now introduce a subgroup of $\Lambda$ for which the symmetry doesn't hold any longer.

\begin{definition}
\label{def - stabilizer}
We denote by ${\rm Stab}_\Lambda(x_\epsilon)$ 
(where $\epsilon \! \in \! \{ +,- \}$)  
the stabilizer in $\Lambda$ of any point $x_\epsilon \!\in\! T_\epsilon$.
We use $\Gamma$ instead of ${\rm Stab}_\Lambda(x_\epsilon)$ when $x_\epsilon$ is the 
vertex $x^-_0$ or $x^-_1$ and the choice of one of these two points is clear.
\end{definition}

The group ${\rm Stab}_\Lambda(x_\epsilon)$ naturally acts on the tree of sign $-\epsilon$.
Since $\Lambda$ acts edge-transitively on $T_-$, up to conjugacy the possible
$\Gamma$-actions are those of
${\rm Stab}_\Lambda(x^-_0)$, of ${\rm Stab}_\Lambda(x^-_1)$ and of their intersection
${\rm Stab}_\Lambda([x^-_0;x^-_1])$.
Henceforth, we fix a vertex $v_-:= x^-_i$, $i \!\in\! \{ 0;1 \}$, and consider the
corresponding $\Gamma$-action on the positive tree $T_+$.
The subgroup generated by the root groups $U_a$ indexed by the twin roots $a$ containing
$v_-$ is a finite index subgroup of $\Gamma=\Lambda(v_-)$.

\subsection{Connection with lattices of Nagao type}
\label{ss - connection with Nagao}
Let us first recall the technical condition on Moufang twin trees already stated in the
introduction.

\begin{definition}
\label{definition - COMM}
We say that a Moufang twin tree $(T_\pm,\delta^*)$, or the associated group $\Lambda$,
satisfies condition $({\rm Comm})$ if any two root groups $U_a$ and $U_b$ commute whenever
$a$ and $b$ are different, prenilpotent twin roots.
\end{definition}

We can now formulate the result relating Moufang twin trees and lattices of Nagao type.


\begin{prop}
\label{prop - Moufang is Nagao}
Let $(T_\pm,\delta^*)$ be a thick locally finite Moufang twin tree.
Let $A$ be its automorphism group and let $\Lambda$ be the subgroup generated by the root
groups.
We fix a twin apartment $L_\pm$, a pair of opposite vertices $v_\pm \!\in\! L_\pm$, and we
denote by $\Gamma$ the group ${\rm Stab}_\Lambda(v_-)$.
For each sign $\epsilon=\pm$, we number
$\{v_j^\epsilon \}_{j \in {\bf Z}}$ the vertices of $L_\epsilon$ in such a way that
$v_0^\epsilon=v_\epsilon$ and for any $j$,
$v_j^\epsilon$ and $v_{j+1}^\epsilon$ are neighbours and $v_j^\epsilon$ and $v_j^{-\epsilon}$
are opposite.
For each $j \!\in\! {\bf Z}$, we denote by $V_j$ (resp. $U_j$) the positive (resp.
negative) root group associated to the twin root whose boundary is $\{v_j^+;v_j^-\}$, and we
denote by $E_+$ (resp. $E_-$) the edge $[v_0^+;v_1^+]$ (resp. $[v_0^-;v_1^-]$).
\begin{enumerate}
\item[(i)] Any of the two geodesic rays $\{v_j \}_{j \geq 0}$ and
$\{v_j \}_{j \leq 0}$ from $v_+$ is a fundamental domain for the
$\Gamma$-action on $T_+$.
\item[(ii)] The stabilizer ${\rm Stab}_A(E_+ \cup \{v_-\})$ is finite, 
and in particular so is $H_0:={\rm Stab}_\Gamma(E_+)$.
Moreover we have the decomposition:
${\rm Stab}_\Gamma(v_j)=H_0 \ltimes (U_1.\cdots.U_j)$ for each $j>0$, and the stabilizer
${\rm Stab}_\Gamma(v_0)$ is the finite group generated by $H_0$ and the root groups indexed
by the two opposite twin roots bounded by $v_-$ and $v_+$.
\item[(iii)] The group $\Gamma$, as well as its action on $T_+$, identifies with the group
and the action attached to the graph of groups of Picture \ref{pict - graph of groups},
where the geodesic ray is $\{v_j \}_{j \geq 0}$ and the groups are the above stabilizers.
\item[(iv)] The group $\Gamma$ is a $T_+$-lattice of Nagao type, which is of directly split
Nagao  type whenever $(T_\pm,\delta^*)$ satisfies $({\rm Comm})$.
\item[(v)] The level function on the positive vertices is nothing else than the codistance
from the $\Gamma$-fixed negative vertex $v_-$, that is: $\ell(x)=\delta^*(v_-,x)$ for any
$x \!\in\! VT_+$.
\end{enumerate}
\end{prop}

{\it Proof.~}
Recall that two edges are called {\it opposite~}if the two corresponding pairs of
vertices of opposite signs and same type are pairs of opposite vertices.
According to \cite[Proposition 7]{Tits92}, the groups $A$ and $\Lambda$ both admit a
structure of RGD-system \cite[3.3]{Tits92}, also called a twin root datum -- see also
\cite[I, Definition 2]{Abr} or \cite[Chapter 1]{RemAst}.
To the RGD-system is naturally attached a positive (resp. negative) Tits system whose Borel
subgroups are the  stabilizers of the edges in $T_+$ (resp. in $T_-$)
\cite[\S 1, Proposition 1]{Abr}.
The group $\Gamma$ is a parabolic subgroup of the negative Tits system and (i) is a
special case of \cite[\S 3, Corollary 1]{Abr}, which applies to general Moufang twin
buildings.

\vskip 3mm

We have $\delta^*(v_1^+,v_0^-)=1$ because $v_1^+$ and $v_1^-$ are opposite and $v_1^-$ and
$v_0^-$ are neighbors.
Since $L_\pm$ is a twin apartment, $v_1^-$ is the only neighbor of $v_0^-$ in $L_-$ which is
opposite $v_1^+$, so that $v_{-1}^-$ is the only neighbor of $v_0^-$ to be at codistance 2
from $v_1^+$.
Therefore the positive root group $V_0$ acts transitively on the neighbors of $v_0^-$ which
are $\neq v_{-1}^-$.
Let $h \!\in\! {\rm Stab}_A(E_+ \cup \{v_-\})$; by the previous sentence, we can find
$u \!\in\! V_0$ such that $u^{-1}h \!\in\! {\rm Stab}_A(E_+ \cup E_-)$, so that
${\rm Stab}_A(E_+ \cup \{v_-\})=V_0 \cdot {\rm Stab}_A(E_+ \cup E_-)$.
Finally, the finiteness of ${\rm Stab}_A(E_+ \cup E_-)$ follows from Ronan-Tits' rigidity
theorem \cite[Theorem 4.1]{RonTit1}, which says that the identity is the only twinning
automorphism which fixes $E_+$, $E_-$ and all the edges having a vertex in common with
$E_+$.
Since $\Lambda < A$, this obviously implies the finiteness of $H_0={\rm Stab}_\Gamma(E_+)$.
The rest of (ii) is a special case of Levi decompositions of stabilizers
of pairs of points of opposite signs in twin buildings \cite[6.3.4]{RemAst}.

\vskip 3mm

For a general Moufang twin tree, the group $U_1.\cdots.U_i$ is in bijection with the set
$U_1 \times \cdots \times U_i$, and it is isomorphic to the direct product group
$U_1 \times \cdots \times U_i$ whenever (Comm) is satisfied: this proves (iv), once we note
that (iii) is a classical consequence of Bass-Serre theory
\cite[I.4.5, Th\'eor\`eme 10]{Ser}.
At last, (v) follows from the fact that the two functions $\ell(\cdot)$ and
$\delta^*(v_-,\cdot)$ coincide on the ray $\{v_j \}_{j \geq 0}$ and are constant on each
$\Gamma$-orbit.
\qed

\begin{rk}\rm
In terms of twin root data (or, equivalently, of RGD-systems), the group $\Delta$ of
\ref{ss - Delta-action on the tree and on horoballs} is the unipotent radical of the
negative parabolic subgroup $\Gamma$, and the decomposition $\Gamma = \Gamma_0 \ltimes \Delta$
is a Levi decomposition \cite[6.2.2]{RemAst}.
The fundamental domain $F$ of Lemma \ref{lemma - fd} is the union of the level-increasing
rays from $v_+$: it is the infinite star ${\rm St}^\infty(v_+)$ defined purely
in terms of codistance in \ref{ss - ttt}.
\end{rk}

Recall that two subgroups of a given group are {\it commensurable~} if they share a finite
index subgroup.

\begin{lemma}
\label{lemma - automorphism in commensurator}
With the same notation as above, we have:
\begin{enumerate}
\item[(i)] For any two points $x, x'\!\in\! T_-$, the groups ${\rm Stab}_A(x)$
and ${\rm Stab}_A(x')$ are commensurable.
\item[(ii)] For any point $x \!\in\! T_-$, we have
$A={\rm Aut}(T_\pm,\delta^*) < C(\Gamma)$.
\item[(iii)] The group $C(\Gamma)$ acts transitively on geodesics of given length and type in
$T_+$, and it acts transitively on geodesics of given length if $T_+$ has a type-exchanging
automorphism.
\end{enumerate}
\end{lemma}

{\it Proof.~}
For any $x \!\in\! T_-$, we can choose an edge $E$ whose closure contains $v$, and in this
case we have: ${\rm Stab}_{A^\circ}(E)<_{\rm f.i.}{\rm Stab}_A(x)$.
Therefore, in order to prove (i) it is enough to consider the case when $x$ and $x'$ both
belong to the interior of an edge, say $E$ for $x$ and $E'$ for $x'$.
Let $y$ (resp. $y'$) be the midpoint of $E$ (resp. $E'$), and let $\delta$ be the distance
between $y$ and $y'$.
Let $g$ be an arbitrary automorphism of $T_+$ stabilizing the edge $E$.
Then $g.y'$ is at distance $\delta$ from $y$ and when $g$ is type-preserving, the geodesic
from $y$ to $g.y'$ leaves $E$ through the vertex through  which the geodesic from $y$ to
$y'$ leaves $E$.
Using root groups for roots containing $E$, we see for any  $g \!\in\! {\rm
Stab}_{G^\circ}(E)$, the edge $g.E'$ is a ${\rm Stab}_\Lambda(E)$-transform of
$E'$.
Since $T_-$ is locally finite, we have a finite set of such edges, say $h_1.E',...$
$h_N.E'$ with $h_1,...$ $h_N \!\in\! {\rm Stab}_\Lambda(E)$.
Therefore ${\rm Stab}_{G^\circ}(E)=\bigsqcup_{j=1}^m h_j
\bigl({\rm Stab}_{G^\circ}(E)\cap{\rm Stab}_{G^\circ}(E')\bigr)$, so that
${\rm Stab}_{A^\circ}(E)\cap{\rm Stab}_{A^\circ}(E')<_{\rm f.i.}{\rm Stab}_{A^\circ}(E)$.
Switching $x$ and $x'$ and using the first remark, we obtain (i).
Moreover (ii) follows from (i) since $g{\rm Stab}_G(x)g^{-1}={\rm Stab}_G(g.x)$ for any
$g \!\in\! G$.
The subgroup $\Lambda$ already enjoys the transitivity property of the first assertion of
(iii) (this is true in the general twin building case \cite[\S 2, p.28]{Abr}), and
the second assertion is easily deduced from the first one.
\qed\vskip 3mm

We can finally turn to the proof of our second main result.

\vskip 3mm

{\it Proof of Theorem \ref{th - dense when Moufang}.~}
We choose $(T_\pm,\delta^*)$ a thick locally finite Moufang twin tree satisfying
{\rm (Comm)}, and go on using the notation of the previous two results.
By Proposition \ref{prop - Moufang is Nagao} (iv), the lattice $\Gamma$ is of
directly split Nagao type.
The closure group $\overline{C(\Gamma)}$ contains $\overline{C(\Gamma) \cap L}$, which is
equal to $L$ by Theorem \ref{th - density}.
But $\overline{C(\Gamma)}$ also contains $\Lambda$, in which many elements do not preserve
the level function $\ell$.
To see this, we note that by Proposition \ref{prop - Moufang is Nagao} (v) we have:
$\ell(x)=\delta^*(v_-,x)$ and $\ell(g.x)=\delta^*(g^{-1}.v_-,x)$ for any
$x \!\in\! VT_+$ and any $g \!\in\! A$.
In order to have $\ell(x) \neq \ell(g.x)$,
it is enough to make $x:=v_+$ and to pick an element $n \!\in\! N_\Lambda(L_\pm)$
lifting a reflection in a vertex $x_j$ with $j \neq 0$ (which always exists e.g. according
to \cite[\S 2, Lemma 4 (ii)]{Abr}).
Finally it remains to apply Proposition \ref{prop - maximality of L} to the group
$F=\overline{C(\Gamma)}$, which contains $L \cup \Lambda$.
\qed

\subsection{Examples}
\label{ss - examples}
Let us now give examples of Moufang twin trees from the most familiar to the most exotic
ones.  We first note that according to M. Ronan and J. Tits \cite[Corollary 8.2]{RonTit2}, a
thick semihomogeneous tree whose set of vertices has cardinality $\alpha$, belongs to
$2^\alpha$ isomorphism classes of twinnings, among which most of them have no automorphism.
The examples below are all required to be Moufang (in particular, have big automorphism
groups).

\begin{example}\rm
There is a detailed study of  $\Lambda={\rm PSL}_2({\bf F}_q[t,t^{-1}])$ in
\cite[\S 2]{RonTit1}.
Let $M$ be a rank 2 free ${\bf F}_q[t,t^{-1}]$-lattice with basis
$\{ u;v \}$, and let $T_+$ and $T_-$ be the Bruhat-Tits trees of
${\rm SL}_2 \bigl( {\bf F}_q(\!(t)\!) \bigr)$ and
${\rm SL}_2 \bigl( {\bf F}_q(\!(t^{-1})\!) \bigr)$,
respectively \cite[II.1]{Ser}.
The vertices in $T_+$, resp. in $T_-$, are the homothety classes of
${\bf F}_q[[t]]$-lattices, resp. of ${\bf F}_q[[t^{-1}]]$-lattices, in
the vector space $M \otimes_{{\bf F}_q[t,t^{-1}]} {\bf F}_q(\!(t)\!)$,
resp. in $M \otimes_{{\bf F}_q[t,t^{-1}]} {\bf F}_q(\!(t^{-1})\!)$.
The codistance is defined in \cite[Lemma 2.1]{RonTit1}.
A useful subset of vertices in $T_+$ is given by the homothety classes of the lattices
${\bf F}_q[[t]] u \oplus {\bf F}_q[[t]] t^j v$ for $j \!\in\! {\bf Z}$, and similarly in
$T_-$ replacing $t$ by $t^{-1}$.
The convex hull of these vertices in $T_+$ is a geodesic $L_+$, and together with the similar
convex hull $L_-$ in $T_-$, they form a twin apartment $L$.
The pointwise fixator $H$ of $L_\pm$ is the finite subgroup of diagonal matrices with
coefficients in ${\bf F}_q^\times$.
The Nagao lattice $\Gamma={\rm PSL}_2({\bf F}_q[t^{-1}])$ lies in $\Lambda$ and the latter
group acts edge-transitively on each tree.
It can actually be proved by hand that for each $d \!\in\! {\bf Z}$, the group
$\{ \left( \begin{array}{rr} 1 & \lambda.t^d \\
0 & 1 \\ \end{array} \right) : \lambda \!\in\! {\bf F}_q \} \simeq ({\bf F}_q,+)$
is the root
group attached to the root $a_d$ defined as the convex hull of the homothety classes of the
lattices ${\bf F}_q[[t]] u \oplus {\bf F}_q[[t]] t^j v$ for $j \geq d$.
The union of these root groups when $d$ ranges over ${\bf Z}$ is the unipotent radical of the
upper triangular Borel subgroup, and condition {\rm (Comm)} holds because this unipotent
group is abelian.
\end{example}

\begin{example}\rm
Split Kac-Moody groups of rank two:
these groups are defined by generators and relations in full generality in \cite{Tit87}.
The datum needed to define such a group roughly consists in a
generalized Cartan matrix $A$ and a groundfield ${\bf K}$.
The rank 2 condition says that the matrix is $2 \times 2$, the size 2 corresponding to the
number of generators of the Weyl group.
The diagonal entries are equal to 2, and one must require that the product of the off-diagonal
entries be $\geq 4$ so that the Weyl group of the building is infinite (dihedral).
Each of the off-diagonal coefficients is negative, and the group satisfies the assumption
{\rm (Comm)} if and only if both entries are $\leq -2$.
This follows from the computation of commutator relations between root groups due, up to
signs, to J. Morita \cite[\S 3, example 6]{Mor}.
Note that apart from the density of the commensurator for the rank 2 case, many other
arguments supporting the analogy between (parabolic subgroups of) Kac-Moody groups over
finite fields and arithmetic groups over function fields are given in \cite{RemNew}.
\end{example}

\begin{example}\rm
Twisted Kac-Moody groups are defined in \cite[Part II, \S 11]{RemAst}.
They are groups consisting of fixed points in a split Kac-Moody group, for a suitable Galois
action.
The rank 2 condition refers to the number of generators in the Weyl group.
A classical example is provided by the unitary group ${\rm SU}_3({\bf F}_q[t,t^{-1}])$, whose
twinned trees are semihomogeneous of valencies $1+q$ and $1+q^3$ \cite[\S 3.5]{RemBie}.
The unipotent radicals of Borel subgroups of ${\rm SU}(3)$ over
${\bf F}_q(\!(t)\!)$ or ${\bf F}_q(\!(t^{-1})\!)$ are not abelian but only metabelian:
${\rm SU}_3({\bf F}_q[t,t^{-1}])$ doesn't satisfy {\rm (Comm)}.
The example of a semihomogeneous Moufang twin tree whose corresponding group $\Lambda$
satisfies {\rm (Comm)~} is given for instance in \cite[13.3]{RemAst}.
The twisted group is a subgroup of rational points in a Kac-Moody group naturally acting on
buildings whose apartments are tilings of the hyperbolic plane.
\end{example}

\begin{example}\rm
\label{ex - not linear but arithmetic}
A general construction of Moufang twin trees is given in \cite[\S 9]{cdf}.
Its starting point are two abstract root group data
of rank 1, and the classification of the latter structure is known in the finite case.
As Tits showed, starting with this data, it is always possible to construct a Moufang
twin tree satisfying our Condition $({\rm Comm})$ and having as root groups (up to
isomorphism) those occurring in the original root group data of rank 1.
An easy way to produce non--classical locally finite Moufang twin trees by this method
was pointed out by Sh. Mozes to the second author.
Taking the affine group ${\bf F}_l^\times \ltimes {\bf F}_l$ ($l$ is any prime power)
acting on the affine line ${\bf F}_l$, we obtain what Tits calls a Moufang set and hence
a root group data of rank 1. We take two copies of this, and then 
J. Tits' method provides homogeneous Moufang twin trees with root
groups isomorphic to ${\bf F}_l^\times$. This cannot occur, at least if $l - 1$ is not 
itself a prime power, for Bruhat-Tits or Kac-Moody
trees since in the latter cases the valencies are cardinalities of projective lines.

\vskip 3mm

Therefore, for instance with $l=7$, we obtain a non-uniform tree
lattice of directly split Nagao type for a regular tree of valency $7$.  In the latter case,
the root groups are all isomorphic to
${\bf F}_7^\times \simeq {\bf Z}/2 \times {\bf Z}/3$.
Let us use the notation of the previous subsection, in particular the $T_+$-lattice $\Gamma$
is the parabolic subgroup of $\Lambda$ fixing a negative vertex $v_-$ of $(T_\pm,\delta^*)$.
The subgroup $H$ of $\Gamma$ generated by the root groups indexed by the negative roots
containing a given geodesic subray of $L_-$ emanating from $v_-$, is isomorphic to the
direct sum of these root groups.
Taking the 2-torsion (resp. 3-torsion) part of $H$ and arguing as in
\cite[Theorem 5]{RemNew}, we conclude that $H$, hence $\Gamma$, cannot be linear over any
field.
In other words, we have obtained a group inclusion $\Gamma < G$ where the
tree-lattice
$\Gamma$ is non-linear but arithmetic in a generalized sense, whereas the classical Nagao
lattice is obviously linear.
Of course, the construction works with any prime power $l$ such that $l-1$ admits two
different prime divisors.
\end{example}

\begin{example}\rm
In J. Tits' construction of the previous example, the starting point is a pair
of two root group data.
In the case where each root group datum comes from a ${\rm PSL}_2$-action
on a projective line, a down-to-earth construction is made in \cite[\S 2]{RemRon}.
This is enough to produce strictly more general Moufang twin trees than those coming from
Kac-Moody groups since there may be two groundfields (one for each type of vertex).
The choices of the characteristics of the ground fields can be made so that the argument of
the second paragraph of the previous example enables to obtain again \og non-linear but
arithmetic\fg tree-lattices.
At last, the concrete viewpoint allows to generalize the construction to the case of Moufang
twin buildings with a Weyl group of arbitrarily large rank, leading to  generalized Kac-Moody
groups with strong non-linearity properties acting on two-dimensional buildings \cite[Theorem
4.A]{RemRon}.
\end{example}

\bibliographystyle{amsalpha}
\bibliography{Trees6.bbl}

\vspace{1cm}

Department of Mathematics \hfill
Institut Fourier -- UMR 5582 du CNRS

University of Virginia \hfill
Universit\'e de Grenoble 1 -- Joseph Fourier

P.O. Box 400137 -- Charlottesville \hfill
100, rue des maths -- BP 74

VA 22904 4137 -- USA \hfill
38402 St Martin d'H\`eres Cedex -- France

{\tt pa8e@virginia.edu} \hfill
{\tt bremy@ujf-grenoble.fr}

\addtolength{\parindent}{-1.6pt}

\end{document}